\documentclass{amsart}
\usepackage{amssymb}
\usepackage{amsthm}
\usepackage{amsmath}
\usepackage{graphicx}
\usepackage{subfigure}
\usepackage{longtable}
\usepackage{psfrag}
\newtheorem{theorem}[subsection]{Theorem}
\theoremstyle{definition}
\newtheorem{definition}[subsection]{Definition}

\newcommand{\F}{\mathbb{F}}
\newcommand{\Z}{\mathbb{Z}}
\newcommand{\Q}{\mathbb{Q}}
\newcommand{\R}{\mathbb{R}}
\newcommand{\C}{\mathbb{C}}
\newcommand{\cH}{\mathcal{H}}
\newcommand{\cV}{\mathcal{V}}
\newcommand{\cW}{\mathcal{W}}
\DeclareMathOperator{\Frob}{Frob}
\DeclareMathOperator{\Gal}{Gal}
\DeclareMathOperator{\im}{im}
\newcommand{\GL}{\mathrm{GL}}
\newcommand{\SL}{\mathrm{SL}}

\begin{document}

\title[Twisted 1-dim cohomology for $\SL(4)$ and Galois representations]
{Cohomology with twisted one-dimensional coefficients for congruence subgroups of $\SL(4,\Z)$ and Galois representations}

\author{Avner Ash} \address{Boston College\\ Chestnut Hill, MA 02445}
\email{Avner.Ash@bc.edu} \author{Paul E. Gunnells}
\address{University of Massachusetts Amherst\\ Amherst, MA 01003}
\email{gunnells@math.umass.edu} \author{Mark McConnell}
\address{Princeton University\\ Princeton, New Jersey 08540}
\email{markwm@princeton.edu}

\thanks{
AA wishes to thank the National Science Foundation for support of this
research through NSF grant DMS-0455240, and also the NSA through
grants H98230-09-1-0050 and H98230-13-1-0261.  This manuscript is
submitted for publication with the understanding that the United
States government is authorized to produce and distribute reprints.
PG wishes to thank the National Science Foundation for support of this
research through the NSF grants DMS-0801214 and DMS-1101640.}



\begin{abstract}
We extend the computations in \cite{AGM1, AGM2, AGM3} to find the
cohomology in degree five of a congruence subgroup~$\Gamma$ of
$\SL(4,\Z)$ with coefficients in a field~$K$, twisted by a nebentype
character~$\eta$, along with the action of the Hecke algebra.  This is
the top cuspidal degree. In practice we take $K=\F$, a finite field of
large characteristic, as a proxy for~$\C$.  For each Hecke eigenclass
found, we produce a Galois representation that appears to be
attached to it.  Our  computations show that in every case this Galois representation is the only one that could be attached to it.
The existence of the attached Galois representations agrees with a theorem of Scholze \cite{scholze}
and sheds light on the Borel-Serre boundary for~$\Gamma$.

The computations require serious modifications to our previous
algorithms to accommodate the twisted coefficients.  Nontrivial
coefficients add a layer of complication to our data structures, and
new possibilites arise that must be taken into account in the Galois
Finder, the code that finds the Galois representations. We have
improved the Galois Finder so that it reports when the attached Galois
representation is uniquely determined by our data.
\end{abstract}

\maketitle

\section{Introduction}\label{intro}

\subsection{} This paper is the next step in our series of papers
\cite{AGM1, AGM2, AGM3} devoted to the computation of the cohomology
of congruence subgroups $\Gamma\subseteq \SL(4,\Z)$ with complex
vector spaces as coefficient modules, together with the action of the
Hecke operators on the cohomology.  In this paper the coefficient
modules are twists of $\C$ by a nebentype character~$\eta$.  The
earlier papers only considered constant coefficients, $\eta = 1$.

We say that a representation $\rho$ of the absolute Galois group of
$\Q$ is attached to a Hecke eigenclass $z$ if for almost all primes
$\ell$, the characteristic polynomial of $\rho(\Frob_\ell)$ is equal
to the Hecke polynomial at $\ell$.  If we only verify the equality
computationally for a finite set of $\ell$, we say that $\rho$
``appears" to be attached to $z$.  In this paper, besides computing
the Hecke operators for small~$\ell$, we find, for each Hecke
eigenclass, a Galois representation that appears to be attached to it.
These Galois representations are uniquely determined, in a sense to be
explained in Section~\ref{introgalreps}.

At the moment, computations of the cohomology of a congruence subgroup
of $\SL(n,\Z)$ as a Hecke module are only feasible when $n$ either
equals the virtual cohomological dimension (vcd) of $\SL(n,\Z)$ or is
one less.  In the latter case we can use Gunnells' algorithm for
computing Hecke operators \cite{experimental}.  When $n=4$, the vcd is
6.  We compute in degree 5 rather than degree 6 because $H^5(\Gamma,
\C)$ supports cuspidal cohomology.

The next paper in this series is planned to deal with
higher-di\-men\-sion\-al twisted coefficients.  This may allow us to
test the generalizations of Serre's conjecture \cite{ADP},
\cite{Herzig} in cases of interest.  The work in this paper is an
important step towards implementing this longer range goal.  In
addition, the current paper already reveals some new phenomena that go
beyond our previous investigations.

In order to avoid the inaccuracy of floating point numbers in our huge
linear algebra computations, we use a finite field~$\F = \F_{p^r}$ as
a proxy for $\C$.  If $p>5$ and there is no $p$-torsion in the
$\Z$-cohomology, then the $\C$- and mod~$p$-betti numbers coincide.
We generally use primes that have four or five decimal digits.
Both~$p$ and the degree~$r$ are chosen to meet certain criteria.  We
choose~$p$ so that the exponent of $(\Z/N)^\times$ divides $p-1$.
This makes the group of characters $(\Z/N)^\times\to\F_p^\times$
isomorphic to the group of characters $(\Z/N)^\times\to\C^\times$.  We
choose~$r$ to ensure that the various Hecke eigenvalues are
$\F$-rational (see Section~\ref{gal}).  Then for any
$\F^\times$-valued character~$\eta$ of $(\Z/N)^\times$, we define
$\F_{\eta}$ to be the one-di\-men\-sion\-al vector space $\F$ regarded
a $\Gamma_0(N)$-module with action via the nebentype character~$\eta$
(Section~\ref{ss:eta}).

Let~$\Gamma$ be a subgroup of finite index in $\SL_4(\Z)$.  To compute cohomology, we use the \emph{sharbly
  complex} $Sh_{\bullet}$, defined in Section~\ref{sh}.  There is an isomorphism of Hecke modules
\[
H^5(\Gamma,M)\approx H_1(\Gamma, Sh_{\bullet}\otimes_{\Z} M),
\]
where $M$ is any module on which the orders of the finite subgroups in
$\SL_{4} (\Z )$ are invertible; this condition is satisfied for us
since we take $M =\F_{\eta }$ to have characteristic $>5$.  In essence, the
method we use to compute $H_1(\Gamma,Sh_{\bullet}\otimes_\Z \F_\eta)$
as a Hecke module is the same as in our papers cited above.  However,
we rewrote the data structures completely to accommodate twisted
coefficients.  The algorithm for the Hecke operators was also
modified.  Since we want in our next papers to compute with twisted
coefficients that are $\F_{p^{r}}$-vector spaces of dimension greater than
one, we made the modifications to accommodate such general coefficient
modules.  The modifications proved to be somewhat tricky.  We explain
them in Sections~\ref{sh} and~\ref{gu}.

\subsection{} \label{introgalreps} Cohomology with twisted
coefficients is interesting because it gives new examples of Scholze's
theorem (recalled in \ref{hlttScholze}) and because we can use it to
test the Serre-type conjectures mentioned above in new cases.  Our
data also raise new questions about the geometry and number theory of
the boundary $B_\Gamma$ of the Borel-Serre compactification $\overline
X/\Gamma$ of the locally symmetric space for~$\Gamma$.  Namely, which
classes in $H^*(\overline X/\Gamma,V)$ restrict to nonzero classes in
$H^*(B_\Gamma,V)$, and why?

To interpret the cohomology properly, we search for Galois
representations attached to Hecke eigenclasses, using a computer
program, the Galois Finder, described in Section~\ref{gal}.  This
program is a modification of the one we used in~\cite{AGM6}.  We look
for reducible Galois representations of degree 4.  For constituents,
we search through all the 1-di\-men\-sion\-al representations
corresponding to Dirichlet characters.  We also search through all the
2-di\-men\-sion\-al representations coming from classical holomorphic
modular forms of weights~2, 3, and~4, as explained in
Section~\ref{ChiReasons}.  Because we look at nontrivial nebentype
characters, classical modular forms of weight 3 can occur.  (With
constant coefficients, only modular forms of weight~ 2 and~4 made an
appearance.)  We also consider symmetric squares of these Galois
representations.  Besides these, for the largest level we consider in
this paper, $N=41$, we needed to employ, as a 3-di\-men\-sion\-al
constituent, a Galois representation which is attached to a cuspidal
automorphic representation with fixed vector under the congruence
subgroup $\Gamma_1(41)$ of $\SL(3,\Z)$ and which is not a symmetric
square lift from $\GL_2$.  (In \cite{AGM2, AGM3} we computed for prime
levels up to $N=211$ and trivial nebentype characters; there we found
other examples of 3-di\-men\-sion\-al constituents that are not
symmetric squares, and examples of irreducible 4-di\-men\-sion\-al
representations attached to Siegel modular forms.)

We always find a unique Galois representation, in the following sense.
For a given Hecke eigenclass $z$, we can compute the Hecke operators
only for a few primes~$\ell$.  (Throughout the paper, $\ell\nmid pN$.)
We always find exactly one Galois representation of the type we search
for that is apparently attached to $z$.  For a more explicit
description of this uniqueness and the procedure used to check it, see
Section~\ref{gal2}.

Our computations are complete for all composite $N\leqslant 28$ and
all prime $N\leqslant 41$.  All the Galois representations found for
levels $N < 41$, whether with trivial or non-trivial~$\eta$, have as
constituents only Dirichlet characters or representations attached to
classical modular forms or their symmetric squares.

It remains unclear why certain combinations of characters and
cusp forms appear in our data and others do not.  We know of no
suffciently explicit computation of the cohomology of the Borel-Serre
boundary for congruence subgroups $\SL(4,\Z)$, and of the Eisenstein
lifting problem for it, that would explain our findings.

The existence of apparently attached Galois representations helps to
corroborate the correctness of our computations.  It is very unlikely
that apparently attached Galois representations could be found if the
computed Hecke eigenvalues were random collections of numbers computed
erroneously.

If we could compute enough Hecke operators, then, using Scholze's
theorem and the method of Faltings-Serre, we could prove that the
apparently attached Galois representations we find are truly attached.
But the computational cost of finding the Hecke eigenvalues at
primes~$\ell$ greater than~$17$ or so is too great, as will be
discussed in the next section.

\subsection{}\label{ss:1.3}

The size of the level $N$ and primes~$\ell$ in the
Hecke operators $T(\ell,k)$ which we compute is limited by
computer speed and memory size. $N$ is limited by the size of the
memory and by the speed, because the numbers of rows and columns of
the matrices which compute the sharbly homology grow like $N^3$.  The
speed limits the Hecke operators, because the number of single cosets
in $T(\ell,k)$ grows like $\ell^3$ for $k=1,3$ and like~$\ell^4$ for
$k=2$.  One new idea in this paper is that, for large~$\ell$, we
compute $T(\ell,1)$ but not $T(\ell,2)$ or $T(\ell,3)$.  This lets us
avoid the~$O(\ell^4)$ part of the computation, while still letting us
eliminate some spurious Galois representations.  See
Section~\ref{gal1}.

\subsection{} \label{hlttScholze}

The work of Peter Scholze \cite{scholze} proves the existence of
attached Galois representations for Hecke eigenclasses in the
cohomology of congruence subgroups of $\GL(n,\Z)$.
This result is conditional on stabilization of the twisted trace
formula.  Scholze's results are easily extended using standard
spectral sequences to coefficient modules which are
finite-di\-men\-sion\-al $\F_p$-vector spaces on which the Hecke
semigroup of matrices acts via reduction modulo some integer $M$.  See
\cite{HLTT} for earlier results of the same kind for
characteristic~0 coefficients.

These results of Scholze allow us to view our cohomology computations
as opening a view onto the world of Galois representations.  Each of
the Hecke eigenclasses we find has an attached Galois representation,
and our computations allow us to investigate exactly which Galois
representations occur, and for which levels and nebentypes.  In this
way we also obtain evidence for the Serre-type conjectures mentioned
above.

\subsection{} Here is a guide to the paper.  In Section~\ref{sh}
we recall the definitions of the Steinberg module, the sharbly
complex, and the concept of attached Galois representation.  In
Section~\ref{gu} we briefly describe how the sharbly homology is
calculated as a Hecke module, with reference to our earlier papers for
details, and with the modifications needed to deal with
$\F_\eta$-coefficients.  In Section~\ref{gal} we describe our Galois
Finder and how it was modified from~\cite{AGM6}.
Section~\ref{ChiReasons} offers an interpretation of our results,
including heuristics.  Section~\ref{res} contains the tables of our
results.

\subsection{Acknowledgments} We thank Darrin Doud, who verified the existence of
the Hecke eigenclass for $\SL_3$ at level $N=41$ that we describe in
Section~\ref{dihedralAndDoud}.  We thank David Rohrlich for helpful
correspondence.

\section{The Sharbly complex, Hecke operators, and Galois representations}\label{sh}

\subsection{}
Let $n\geqslant 2$.  Let $\Q^n$ denote the space of
$n$-di\-men\-sion\-al column
vectors.

\begin{definition} \label{defSharb}
The \emph{Sharbly complex} $Sh_{\bullet } $ is the complex of left $\Z
\GL(n,\Q)$-modules defined as follows.  As an abelian group, $Sh_{k}$
is generated by symbols $[v_1,\dots,v_{n+k}]$, where the $v_i$ are
nonzero vectors in $\Q^n$, modulo the submodule generated by the
following relations:

(i) $[v_{\sigma
(1)},\dots,v_{\sigma(n+k)}]-(-1)^\sigma[v_1,\dots,v_{n+k}]$ for all
permutations $\sigma$;

(ii) $[v_1,\dots,v_{n+k}]$ if $v_1,\dots,v_{n+k}$ do not span all
of $\Q^n$; and

(iii) $[v_1,\dots,v_{n+k}]-[av_1,v_{2},\dots,v_{n+k}]$ for all $a\in
\Q^\times$.

\noindent The element $g \in \GL(n,\Q)$ acts on $Sh_{\bullet}$ by
$g [v_1,\dots,v_{n+k}] = [g v_1,\dots,g v_{n+k}]$.  The boundary map
$\partial_k \colon Sh_{k} \rightarrow Sh_{k-1}$ is
\[ 
\partial_k([v_1,\dots,v_{n+k}])=
\sum_{i=1}^{n+k} (-1)^i [v_1,\dots,\widehat{v_i},\dots v_{n+k}],
\]
where as usual $\widehat{v_i}$ means to delete $v_{i}$.
\end{definition}

All these objects depend on $n$, which we suppress from the
notation, since we will later work only with $n=4$.

The sharbly complex 
$$
\dots\to Sh_i\to Sh_{i-1}\to \dots \to Sh_1\to Sh_0
$$ is an exact sequence of $\GL(n,\Q)$-modules.  We may define the
Steinberg module $St$ as the cokernel of $\partial_1\colon Sh_1\to
Sh_0$ (cf.~\cite[Theorem 5]{AGM5}).

Let $\Gamma$ be a congruence subgroup of $\SL(n,\Z)$.

\begin{definition} Let $M$ be a left $\Gamma$-module.  The
\emph{sharbly homology} of $\Gamma$ with coefficients in $M$ is
$H_*(\Gamma,Sh_{\bullet} \otimes_\Z M)$, where $\Gamma$ acts
diagonally on the tensor product.
\end{definition}

If $(\Gamma,S)$ is a Hecke pair in $\GL(n,\Q)$ and $M$ is a left
$S$-module, the Hecke algebra $\cH(\Gamma,S)$ acts on the sharbly
homology, since $S$ acts (diagonally) on $Sh_{\bullet} \otimes_\Z M$
and because the sharbly homology is the homology of the complex
$H_0(\Gamma, Sh_\bullet \otimes_Z M)$.

The following theorem is proved in~\cite{AGM4}.

\begin{theorem}\label{S-Steinberg}
For any $\Gamma\subset \GL(n,\Z)$ and any coefficient module $M$ in
which all the torsion primes of $\Gamma$ are invertible, there is a
natural isomorphism of Hecke modules
$$
H_i(\Gamma,Sh_{\bullet} \otimes_\Z M)\to  H^{\binom{n}{2}-i}(\Gamma, M)
$$
for all~$i$.
\end{theorem}

\subsection{} \label{ss:eta} We now define the~$\Gamma$ and $\Gamma$-modules used in this paper.

\begin{definition} Let $\Gamma_0(N)$ be the subgroup of matrices in $\SL(n,\Z)$ whose
bottom row is congruent to $(0,\dots,0,{*})$ modulo~$N$.
\end{definition}

Let $\F = \F_{p^r}$ be a finite field of characteristic $p$.  Let
$\eta:(\Z/N)^\times\to\F^\times$ be a character, which we will call
the \emph{nebentype} (even if it is trivial, although in that case we
will sometimes speak of the ``trivial character''.)  In practice, $p$
will be a prime of four or five decimal digits.  We will always
choose~$p$ so that the exponent of $(\Z/N)^\times$ divides $p-1$.
Hence~$\eta$ takes values in~$\F_p^\times$.

Define $S_{pN}$ to be the subsemigroup of integral matrices in
$\GL(n,\Q)$ satisfying the same congruence condition as $\Gamma_0(N)$
and having positive determinant relatively prime to~$pN$.  Let
$\cH(pN)$, the \emph{anemic Hecke algebra}, be the $\Z$-algebra of
double cosets $\Gamma_0(N)S_{pN}\Gamma_0(N)$.  Then $\cH(pN)$ is a
commutative algebra that acts on the cohomology and homology of
$\Gamma_0(N)$ with coefficients in any $S_{pN}$-module.  In
particular, $\cH(pN)$ contains all double cosets of the form
$\Gamma_0(N)D(\ell,k)\Gamma_0(N)$, where $\ell$ is a prime not
dividing $pN$, $0\leqslant k\leqslant n$, and $D(\ell,k)$
is the diagonal matrix with the first $n-k$ diagonal entries equal to
1 and the last $k$ diagonal entries equal to $\ell$.  These double
cosets generate $\cH (pN)$ (cf.~\cite[Thm.~3.20]{Shimura}).  When we
consider the double coset generated by $D(\ell,k)$ as a Hecke
operator, we call it $T(\ell,k)$.

Write ${\F_\eta}$ for the $S_{pN}$-module where a matrix $s\in S_{pN}$
acts on $\F$ via $\eta(s_{nn})$, where~$s_{nn}$ is the~${*}$ in the
bottom row congruent to $(0, \dots, 0, {*})$ mod~$N$.

\begin{definition}\label{def:hp}
Let $V$ be an $\F[\cH(pN)]$-module.  Suppose that $v\in V$ is a
simultaneous eigenvector for all $T(\ell,k)$ and that
$T(\ell,k)v=a(\ell,k)v$ with $a(\ell,k)\in\F$ for all prime $\ell
\nmid pN$ and $0\leqslant k\leqslant n$.  If
\[
\rho\colon G_\Q\to \GL(n,\F)
\]
is a continuous representation of
$G_{\Q} = \Gal (\overline\Q/\Q)$ unramified outside $pN$, and if
\begin{equation}\label{eqn:hp}
\sum_{k=0}^{n}(-1)^k\ell^{k(k-1)/2}a(\ell,k)X^k=\det(I-\rho(\Frob_\ell)X)
\end{equation}
for all $\ell\nmid \ pN$, then we say that~$\rho$ is \emph{attached} to~$v$. 
\end{definition}  
Here, $\Frob_\ell$ refers to an arithmetic Frobenius element, so that
if $\varepsilon$ is the cyclotomic character, we have $\varepsilon
(\Frob_\ell)=\ell$.  

The polynomial in~(\ref{eqn:hp}) is called the \emph{Hecke polynomial}
for~$v$ and~$\ell$.

As explained in the introduction, we have the following special case
of a theorem of Scholze (conditional on stabilization of the twisted
trace formula):

\begin{theorem}\label{scholze} 
Let $N\geqslant1$.  
Let $v$ be a Hecke eigenclass in
$H^{5}(\Gamma_0(N), \F_\eta)$.  Then
there is attached to~$v$ a continuous representation unramified
outside $pN$:
\[
\rho\colon G_\Q\to \GL(n, \F).
\]
\end{theorem}

Echoing Definition~\ref{def:hp}, we say that~$\rho$ is
\emph{apparently attached} to~$v$ if condition~(\ref{eqn:hp}) holds
for a finite range of~$\ell$ which we have computed, a range large
enough that we are confident~$\rho$ really is attached to~$v$.

\section{Computing homology and the Hecke action}\label{gu}

Following Theorem 2.4, we compute the Hecke operators acting on
sharbly cycles that are supported on Voronoi sharblies.  Theorem 13 of
\cite{AGM5} guarantees that the packets of Hecke eigenvalues we compute
do occur on eigenclasses in $H_{1}(\Gamma_0(N),Sh_{\bullet} \otimes_\Z
\F_\eta) \approx H^5(\Gamma,\F_\eta)$.  In this section, we define the
Voronoi sharblies, recall results from \cite{AGM4,AGM5}, and explain
how the results are modified to work with $\F_\eta$-coefficients.

\subsection{}
The sharbly complex is not finitely generated as a $\Z
\SL(n,\Z)$-module, which makes it difficult to use in practice to
compute homology.  To get a finite complex to compute $H_{1}$, we use
the Voronoi complex.  We refer to \cite[Section 5]{AGM5} for any
unexplained notation in what follows.

Let $X_{n}^{0}\subset \R^{n(n+1)/2}$ be the convex cone of
positive-definite real quadratic forms in $n$ variables.  This has a
partial (Satake) compactification $(X_{n}^{0})^{*}$ obtained by
adjoining rational boundary components, and the compactification is
itself a convex cone.  The space $(X_{n}^{0})^{*}$ can be partitioned
into cones $\sigma =\sigma (x_{1},\dotsc ,x_{m})$, called
\emph{Voronoi cones}, where the $x_{i}$ are contained in certain
subsets of nonzero vectors from $\Z^{n}$.  (We write elements of
$\Z^{n}$ as column vectors, as we did in Section~\ref{sh} for
$\Q^{n}$.)  The cones are built as follows.  Each nonzero $x_{i}\in
\Z^{n}$ determines a rank-one quadratic form $q (x_{i}) = x_{i}
x_{i}^t \in (X_{n}^{0})^{*}$.  Let $\Pi$ be the closed convex hull of
the points $\{ q (x)\mid x\in \Z^{n}, x\ne 0\}$.  Then each of the
proper faces of $\Pi$ is a polytope, and the~$\sigma$ are the cones on
these polytopes.  The indexing sets are constructed in the obvious
way: if $\sigma$ is the cone on $F\subset \Pi$, and $F$ has distinct
vertices $q (x_{1}),\dotsc ,q (x_{m})$, then the indexing set is
$\{\pm x_{1},\dotsc ,\pm x_{m} \}$.  We let $\Sigma$ denote the set of
all Voronoi cones.

Let $X^{*}_{n}$ be the quotient of $(X^{0}_{n})^{*}$ by homotheties.
The images of the Voronoi cones are cells in $X^{*}_{n}$.  Let $\Z
V_{\bullet}$ be the oriented chain complex on these cells, graded by
dimension.  Let $\Z \partial V_{\bullet}$ be the subcomplex generated
by those cells that do not meet the interior of $X^{*}_{n}$ (i.e., do
not meet the image in $X^{*}_{n}$ of the positive-definite cone).  The
\emph{Voronoi complex} is then defined to be $\cV_{\bullet} = \Z
V_{\bullet}/\Z \partial V_{\bullet}$.  For our purposes, it is
convenient to reindex $\cV_{\bullet}$ by introducing the complex
$\cW_{\bullet}$, where $\cW_{k} = \cV_{n+k-1}$.  The results of
\cite{AGM4,AGM5} show that, if $n\leqslant 4$, both $\cW_{\bullet}$
and $Sh_{\bullet}$ give resolutions of the Steinberg module.  In
particular, let $\Gamma =\Gamma_{0} (N)$.  If $M$ is a
$\Z[\Gamma]$-module such that the order of all torsion elements in
$\Gamma$ is invertible, then $H_{*} (\Gamma , \cW_{\bullet }
\otimes_{\Z}M)\approx H_{*} (\Gamma , Sh_{\bullet}\otimes_{\Z}M)$, and
furthermore by Borel--Serre duality these are isomorphic (after
reindexing) to $H^{*} (\Gamma , M)$.  These two complexes can be
related as follows when $n=4$:  every Voronoi
cell in $X_4^*$ of dimension $\leqslant 5$ is a simplex.  Thus for $0\leqslant k
\leqslant 2$, we can define a map of $\Z[\SL(4,\Z)]$-modules
$$
\theta_k\colon \cW_{k} \to Sh_{k}
$$ that takes the Voronoi cell $\sigma(v_1,\dots,v_{k+4})$ to
$\theta_k((v_1,\dots,v_{k+4})):=[v_1,\dots,v_{k+4}]$.  This allows us
to realize Voronoi cycles in these degrees in the sharbly complex.
The image of~$\theta_k$ is the set of \emph{Voronoi sharblies} in
degree~$k$.  Then $H_{1} (\Gamma , \cW_{\bullet }
\otimes_{\Z}M)\approx H_{1} (\Gamma , Sh_{\bullet}\otimes_{\Z}M)$ by
\cite[Corollary~12]{AGM5}.

\subsection{} \label{mmCode} We now explain concretely
how we compute $H_{1}(\Gamma,\cW_{\bullet} \otimes_\Z \F_\eta)$.  We
have a body of code in Sage \cite{sage} for these computations.  The
code supports $G$-modules~$M$, that is, representations of~$G$.  Here
$G$ is a finite group, or a matrix group like $\Gamma_0(N)$ or
$S_{pN}$. The module~$M$ has finite dimension over its base ring.  The
base ring is $\F$, $\Q$, or~$\Z$ in this project, though it could be
more general. Morphisms of $G$-modules are supported, as are kernel,
cokernel, image, direct sum, and tensor products of $G$-modules.
When~$H$ is a finite-index subgroup of~$G$, we support
$\mathrm{Res}_H^G$, $\mathrm{Ind}_H^G$, and $\mathrm{Coind}_H^G$ of
$G$-modules, functorially.

The program takes as input the values of~$N$, $p$, and the nebentype,
which is a one-di\-men\-sion\-al representation~$\eta$ of $S_{pN}$ with
coefficients in~$\F_p$.  (The extension from $\F_p$ to~$\F_{p^r}$
comes later, in the Galois Finder.)  The nebentype is essentially a
Dirichlet character $(\Z/N)^\times \to \F_p^\times$.  Sage makes it
automatic to enumerate the Dirichlet characters.

The complex $\cW_\bullet$ has only finitely many classes of Voronoi
cells modulo $\SL(n,\Z)$ \cite{Vor}.  When $n=4$, there are~$18$
classes.  In fact, to compute $H_1$ we only need $\cW_0$, $\cW_1$, and
$\cW_2$, so our code truncates away the rest of $\cW_\bullet$ for
efficiency.

For each class of cells modulo $\SL(n,\Z)$, the code maintains a
standard representative cell~$\sigma$ as listed in \cite{mmThesis}.
The stabilizer $G_\sigma$ of~$\sigma$ in $\SL(n,\Z)$ acts on~$\sigma$
with orientation character $Z_\sigma$.  The code stores~$G_\sigma$
and~$Z_\sigma$.

Fix right coset representatives~$r$, $r'$, \dots\ for
$\Gamma_0(N)\backslash\SL(n,\Z)$ once and for all.  Since
$\Gamma_0(N)$ has finite index in $\SL(n,\Z)$, the complex
$\cW_\bullet$ has only finitely many classes of cells modulo
$\Gamma_0(N)$.  For each class modulo $\Gamma_0(N)$, we may choose a
representative cell $\sigma_1 = r\sigma$, where $\sigma$ is one of the
representative cells modulo $\SL(n,\Z)$, and~$r$ is one of the
standard coset representatives.  An awkward fact is that, for two
different coset representatives $r$, $r'$, the cells $r\sigma$ and
$r'\sigma$ may be in the same $\Gamma_0(N)$-orbit.  This
occurs when $r^{-1} r'$ is in the stabilizer~$G_\sigma \subset
\SL(n,\Z)$.  For computation we must choose $r$ or $r'$, not both; say
we choose~$r$.  A class \texttt{CellOrbitStructure} in our code
handles these details.  $\sigma_1$ itself may have a non-trivial
stabilizer $G_{\sigma_1} \subset \Gamma_0(N)$; the
\texttt{CellOrbitStructure} takes care of these stabilizers
$G_{\sigma_1}$ and how their orientation characters $Z_{\sigma_1}$
interact with the orientation characters $Z_\sigma$ of $G_\sigma$.

Equation~(\ref{celloProb}) below presents a problem we need to
solve repeatedly during the homology calculation.  Suppose we are
given a cell $\tau \in \cW_\bullet$, with $\tau = g \sigma$
for~$\sigma$ a standard cell and for $g\in\SL(n,\Z)$.  Then $g =
\gamma_1 r'$ for some coset representative~$r'$ and $\gamma_1 \in
\Gamma_0(N)$.  Since we chose~$r$ instead of~$r'$, we have $\gamma r'
= r g_\sigma$ for some stabilizer element $g_\sigma \in G_\sigma$ and
some $\gamma\in\Gamma_0(N)$.  Thus
\begin{equation} \label{celloProb}
g = (\gamma_1 \gamma^{-1}) r g_\sigma.
\end{equation}
The problem is, given $g$ and $\sigma$, to solve for $\gamma_1$,
$\gamma$, $r$, $g_\sigma$, and to compute the orientation characters.
The \texttt{CellOrbitStructure} has a method \texttt{decompose} that
solves~(\ref{celloProb}).

Let~$\sigma_1$ run through all the representatives $r\sigma$ of the
classes of cells modulo $\Gamma_0(N)$.  During the homology
computation, we need, for each~$\sigma_1$, to restrict the
nebentype~$\eta$ to the finite stabilizer group $G_{\sigma_1}$, and to
tensor the restriction with the orientation character of
$G_{\sigma_1}$.  This tensor product $\eta_{\sigma_1} : G_{\sigma_1}
\to \F_p^\times$ is called the \emph{local representation}
for~$\sigma_1$.  The \texttt{CellOrbitStructure} keeps track of the
local representations.

As we explain in \cite{AGM4,AGM5}, $H_\bullet (\Gamma, \cW_k
\otimes_{\Z}\F_\eta)$ is computed by a spectral sequence.  The columns
are indexed by~$k$, and the $j$-th row is the direct sum of the
homology groups $H_j(G_{\sigma_1}, \eta_{\sigma_1})$.  Since the
torsion in $\Gamma_0(N)$ has order prime to~$p$, all the homology
groups $H_j$ vanish for $j > 0$.  The $E^1$ term has only one row,
whose entry in the $k$-th box is the module of co-invariants
\[
E^1_{k,0} = H_0(\Gamma, \cW_k \otimes_{\Z} \F_\eta) = \cW_k
\otimes_{\Z\Gamma} \F_\eta.
\]
As~$\sigma_1$ runs through representatives of the cells modulo
$\Gamma_0(N)$, the co-invarant module $\cW_k \otimes_{\Z\Gamma}
\F_\eta$ breaks up as a direct sum:
\begin{equation}\label{celloDel1}
E^1_{k,0} = \bigoplus_{\sigma_1\mathrm{\ of\ degree\ }k} H_0(G_{\sigma_1}, \eta_{\sigma_1}).
\end{equation}
Each summand $H_0(G_{\sigma_1}, \eta_{\sigma_1})$ is the module of
co-invariants for the local representation $\eta_{\sigma_1}$.  It is
isomorphic to~$\F$ if $\eta_{\sigma_1}$ is a trivial representation,
and is zero otherwise.

The $E^2_{k,0}$ of the spectral sequence is isomorphic to $H_k
(\cW_{\bullet } \otimes_{\Z\Gamma}\F_\eta)$.  This is computed using
the differential~$\bar\partial_k$ that is the tensor product
with~$\eta$ of the differential~$\partial_k$ on sharblies in
Section~\ref{defSharb}.  $\bar\partial_k$ is constructed in Sage as a
sparse matrix of size $\dim E^1_{k,0} \times \dim E^1_{k-1,0}$.  As
before, we are computing~$H_1$, so we only compute~$\bar\partial_2$
and~$\bar\partial_1$.

We illustrate the sizes of these matrices with the example of $N=41$,
$p=21881$, and trivial nebentype.  Here $\bar\partial_2$ is
$24590\times 7100$, and $\bar\partial_1$ is $7100 \times 746$.  (This
is small compared to \cite{AGM3}, where, for $N=211$ and trivial
nebentype, $\bar\partial_2$ was about four million by one million.  We
did not compute the Hecke operators in \cite{AGM3}.)

We write the matrices~$\bar\partial_2$ and~$\bar\partial_1$ to disk,
partly as insurance in case of a computer crash during a long run.
The next step is to choose a basis~$\{x_i\}$ of the homology,
$\ker(\bar\partial_1) / \im(\bar\partial_2)$.  We choose the basis
using Sheafhom, a package written by one of us (MM) in Common Lisp and
described in~\cite{AGM3}.  Sheafhom performs homology calculations by
row- and column-reducing large sparse matrices while saving the
change-of-basis matrices to disk.  It works with base rings~$\F_p$ as
well as~$\Z$.  If~$y$ is a cycle in the homology, Sheafhom can express
it as a linear combination of the homology basis, $y = \sum c_i x_i$,
using only a small amount of RAM.

\subsection{} \label{pgCode} To compute the Hecke operators, we use the basis
$\{x_{i} \}$ we found for the homology group $H_{1} (\cW_{\bullet }
\otimes_{\Z\Gamma }\F_\eta)$.  We identify the~$x_{i}$ with elements
$y_{i} = \theta_{1,*} (x_{i})\in H_{1} (Sh_{\bullet
}\otimes_{\Z\Gamma}\F_\eta)$.  Let $T$ be a Hecke operator.  Using the algorithm of Gunnells mentioned in the introduction, we
compute each Hecke translate $Ty_{i}$ and then find a sharbly cycle
$z_{i}$ such that $z_{i}=Ty_{i}$ in $H_{1}
(Sh_{\bullet}\otimes_{\Z\Gamma}\F_\eta)$ and such that $z_{i}$ is in
the image of the map $\theta_{1,*}$.  The inverse images
$\theta_{1,*}^{-1} (z_{i})$ can be written as linear combinations
$\sum c_i x_i$ as in the previous paragraph.  This gives a matrix
representing the action of $T$.  From this matrix we can find
eigenclasses and eigenvalues.

\section{Finding attached Galois representations}\label{gal}

From now on we set $n=4$.

Suppose we have a finite-di\-men\-sion\-al $\F$-vector space~$V$
together with an action of the Hecke operators from $\cH(pN)$.  We now
describe how we find Galois representations that are apparently
attached to Hecke eigenvectors $v$ in $V$.  Our Galois Finder program
is part of our Sage code.

\subsection{} \label{gal1} As in Section~\ref{pgCode}, we compute 
the action on $V=H_1(\Gamma_0(N), \cW_{\bullet} \otimes_\Z \F_\eta)$
of the Hecke operators $T(\ell,k)$ for $k=1, 2, 3$ and for~$\ell$
ranging through a set
\[
L = \{\ell \mid \ell \mathrm{\ prime,} \, \ell \leqslant \ell_0, \,
\ell\nmid pN\}.
\]
The upper bound~$\ell_0$ depends on the level~$N$ and the
nebentype~$\eta$, because sometimes we need more~$\ell$ to find a
unique Galois representation.  $\ell_0$ is never less than~$5$, and is
occasionally as high as~$17$.  For the larger $\ell$, as we have
mentioned, we sometimes compute only $T(\ell,1)$ and not $T(\ell,k)$
for $k=2,3$.  The tables in Section~\ref{resMain} list which operators
we computed.  $T(\ell,0)$ is always the identity and $T(\ell,4)$ is
$\eta(\ell)$ times the identity.  To check the work, we always verify
that our Hecke operators commute pairwise.

For a given level~$N$, we look for reducible Galois representations
apparently attached to a given packet of Hecke eigenvalues.  For some
$N$ there are irreducible Galois representations attached to certain
packets, but our data does not extend to such large $N$.

Some constituents of the Galois representations we are looking for are
1-di\-men\-sion\-al, coming from Dirichlet characters mod~$N$ taking
values in the cyclotomic field $K_0 = \Q(\zeta_N)$, with
$\zeta_N\in\C$ a primitive $N$-th root of unity.  Others are
2-di\-men\-sion\-al, coming from newforms of level dividing~$N$ and
weights 2, 3, or~4.  Using Sage, we compute the newforms of
level~$N_1$ and character~$\psi$, for all $N_1 \mid N$, all
Dirichlet characters~$\psi$ mod~$N_1$, and all weights 2, 3, and~4
(see Section~\ref{gal3} for full details).  Other constituents are
3-di\-men\-sion\-al, coming either from symmetric squares of
2-di\-men\-sion\-al representations or from $\GL(3)$-homology classes
which are not symmetric squares (which in our data occurs only for
$N=41$).

Consider the fields of definition $K_1, K_2, \dots$ of the newforms we
have listed, together with~$K_0$.  The Galois Finder will be
computing, not in $\F_p$, but in the residue class fields for the
primes $\mathfrak{P}$ over~$p$ in the different~$K_i$.  We define~$r$
to be the smallest integer so that all these residue class fields
embed in $\F = \F_{p^r}$.  The field~$\F$ is recorded at the top of
each table as $\F = GF(p^r)$.

Computation in $\F_{p^r}$ slows down when~$r$ becomes large.  We would
have liked to choose~$p$ so that it splits completely in all
the~$K_i$, meaning $r=1$.  As we have mentioned, we always choose~$p$
so that it splits completely in~$K_0$.  But for $N$ in the 20s and
higher, the fields $K_1$, $K_2$ become large enough that we cannot
choose a four or five digit~$p$ that splits completely everywhere.
Instead, we choose~$p$ so that~$r$ will be as small as possible.

The Hecke operators $T(\ell,k)$ we compute are all semisimple.  We do
not know how to prove that this would always be the case.

For each $T(\ell,k)$ that we compute, we
decompose~$V$ into eigen\-spaces under that operator.  In principle,
the eigen\-values $a(\ell, k)$ of $T(\ell,k)$ might lie in an
extension field of $\F$, but we always observe that they lie in~$\F$.

After decomposing~$V$ into eigen\-spaces, we take the common
refinement of the decompositions.  Let $E$ have the form
$\bigcap_{(\ell, k)} E_{\ell, k}$, where $E_{\ell, k}$ is any one of
the eigenspaces for the operator at $(\ell, k)$, and the intersection
is over all $\ell\in L$ and all~$k$ we have computed.  We find all the
non-zero~$E$ of this form.  They are the \emph{simultaneous
  eigen\-spaces}.  $V$ is the direct sum of the~$E$.  By construction,
the Hecke eigenvalues $a(\ell, k)$ are constant on each~$E$ and
characterize it.  The function $(\ell, k) \mapsto a(\ell, k)$ is the
\emph{Hecke eigenpacket} of~$E$.

We distinguish two kinds of multiplicity for~$E$.  We define the
\emph{Hecke multiplicity} of~$E$ to be $\dim_{\F} E$.  A second kind,
the \emph{Galois multiplicity}, is defined in Section~\ref{galmult}.

To a simultaneous eigenspace~$E$ we attach a family of polynomials.
The \emph{polynomial system} $\mathcal{F}(E)$ is the mapping that
sends $\ell \in L$ to the Hecke polynomial with eigenvalues $a(\ell,
k)$ defined in~(\ref{eqn:hp}), or to a partial Hecke polynomial which
we now explain.  For small~$\ell$, we can compute the Hecke
eigenvalues $a(\ell,k)$ for all $k=0,\dots,4$, so we know the whole
Hecke polynomial~(\ref{eqn:hp}); call this a \emph{full} Hecke
polynomial.  For larger~$\ell$, computing $T(\ell,2)$ would be too
slow.  In this case, we compute $T(\ell,1)$, and we only know that the
Hecke polynomial is $1 - a(\ell,1)X + O(X^2)$, where $O(X^2)$ means some undetermined  linear combination of $X^2, X^3,$ and $X^4$.  We call the latter a
\emph{partial} Hecke polynomial.  A partial Hecke polynomial is
implemented in Sage as an element of the quotient ring $\F[X]/(X^2)$.
As a whole, $\mathcal{F}(E)$ contains one or more full polynomials,
all of degree~4, and zero or more partial polynomials, whose degree is
undefined.  We say $\deg\mathcal{F}(E) = 4$.

\subsection{} \label{gal2} We use known Galois representations~$\rho$
unramified outside~$pN$, taking values in $\GL(m, \F)$ for $m = 1$
or~2.  These are the Galois representations coming from Dirichlet
characters and newforms as described roughly in Section~\ref{gal1},
and to be described in full detail in Section~\ref{gal3}.  We also use
the symmetric squares of the~$\rho$ coming from newforms; these take
values in $\GL(m, \F)$ for $m=3$.  The characteristic polynomial of
Frobenius for each of these representations is known and is of
degree~$m$ for each $\ell \nmid pN$.  In the language above, they are
all \emph{full} polynomials.  In our code, we define the
\emph{polynomial system} $\mathcal{F}(\rho)$ to be the mapping that
sends $\ell \in L$ to the characteristic polynomial of Frobenius
for~$\rho$ at~$\ell$.  We say $\deg\mathcal{F}(\rho) = m$.

Let us describe how we conjecturally attach a sum of $\rho$'s to a
simultaneous eigenspace~$E$.  On polynomial systems, it is natural to
define $\mathcal{F}(\rho_1\oplus\cdots\oplus\rho_t) = \prod_{i=1}^t
\mathcal{F}(\rho_i)$, a product of polynomial systems.  We can also
define quotients, but we must be careful about the partial Hecke
polynomials, as we now explain.  Let $\mathcal{F}_1$ and
$\mathcal{F}_2$ be two polynomial systems with the same~$L$.  Say that
$\mathcal{F}_1$ \emph{divides} $\mathcal{F}_2$ if, for each $\ell\in
L$, the polynomial at~$\ell$ for $\mathcal{F}_1$ divides the
polynomial at~$\ell$ for~$\mathcal{F}_2$.  Implicit in this definition
is that $\deg\mathcal{F}_1 \leqslant \deg\mathcal{F}_2$.  When one
polynomial system divides another, define the \emph{quotient system}
in the obvious way.  The degree of the quotient system is
$\deg\mathcal{F}_2 - \deg\mathcal{F}_1$.  For some~$\ell$ we will be
dividing a partial Hecke polynomial by a full Hecke polynomial, but we
never use a partial polynomial as a divisor.  Dividing a partial Hecke
polynomial $f_1(x) \pmod {X^2}$ by a full Hecke polynomial $f_2(x)$ is
well defined because $f_2(x)$ always has constant term~1,
hence~$f_2(x)$ projects via $\F[X] \to \F[X]/(X^2)$ to a unit
$\F[X](X^2)$; the inverse of $1 - aX$ in $\F[X]/(X^2)$ is $1 + aX$.
Note that we can divide indefinitely in $\F[X](X^2)$, because
$1/(1-aX)^\nu$ exists for arbitrarily large~$\nu$.  The reason we keep
track of the degree of a polynomial system is that, although we could
divide into the partial polynomials indefinitely, we will stop
dividing by~$\mathcal{F}_1$ as soon as the full polynomials of the
quotient reach degree~0.

For a given~$E$, we make a list~$\mathcal{R}$ of all the~$\rho$ for
which $\mathcal{F}(\rho)$ divides $\mathcal{F}(E)$.  Then we run
through all possible finite subsets of~$\mathcal{R}$, say $\{\rho_1,
\dots, \rho_t\}$, and we make a list $\mathcal{R}'$ of all the direct
sums $\rho_1\oplus\cdots\oplus\rho_t$ for which
$\mathcal{F}(\rho_1\oplus\cdots\oplus\rho_t) = \mathcal{F}(E)$.  We
always find that $\mathcal{R}'$ is non-empty.  If $\mathcal{R}'$ has
two or more elements, we take more primes~$\ell$, add them to $L$,
compute the Hecke operators $T(\ell,k)$ (or at least $T(\ell,1)$), and
refine the eigenspaces~$E$ for the new operators if necessary.  We
almost always find we can take enough~$\ell$ to make $\mathcal{R}'$
have exactly one element.  In the ``minor" exceptions to this
statement, recounted in the next paragraph, it is still true that the
Galois representations $\rho_1\oplus\cdots\oplus\rho_t$ for the
elements in $\mathcal{R}'$ are isomorphic to each other.  Therefore in
every case we can discover the uniqueness of the Galois representation
among those our Finder looks through that seems to be attached to any
given Hecke eigenspace we have computed.  We assert this uniqueness
even though our data is rather limited, i.e., $L$ is not that large.
Of course, by Chebotarev Density the truly attached Galois
representation is unique, up to semisimplification.

\subsection{} \label{dihedralAndDoud} There are some 
exceptions to the statement that $\mathcal{R}'$ has exactly one
element.  The minor exceptions are found in the tables in
Section~\ref{resMain} at
\begin{itemize}
\item $N=24$, $\eta = \chi_{24,0}\chi_{24,1}\chi_{24,2}$,
  representations with $\sigma_{24,2c}$;
\item $N=28$, $\eta = \chi_{28,0}\chi_{28,1}^3$, representations with
  $\sigma_{28,2c}$.
\end{itemize}
Here, one Galois representation with a symmetric square in it happens
to coincide with one representation without a symmetric square.  We
checked by computer that the Hecke polynomials match for all $\ell <
1000$, $\ell\nmid N$.  In these two cases, the symmetric square is of
a ``dihedral'' Galois representation, so that its symmetric square is
reducible.  These are in fact the same four-di\-men\-sion\-al Galois
representation.

The major exceptions occured at level $N=41$ and the nebentype~$\eta =
\chi_{41}^{10}$ whose image has order~$4$.  Here $\dim V = 8$,
splitting into eight $E$'s of dimension~1, and $\mathcal{R}'$ was
empty for two out of the eight~$E$.  Darrin Doud, upon our request,
using computer programs he developed, found an autochthonous form for
$\SL_3$.  Specifically, he found a three-di\-men\-sion\-al Galois
representation~$\delta$ attached to a cohomology class $z$ for a
congruence subgroup of $\SL(3,\Z)$ and with coefficients in $\F_\eta$,
which is not a lift from any lower-rank group.  Of the two
four-di\-men\-sion\-al representations that we could not identify
using characters and cusp forms, one proved to be $1 \oplus
\varepsilon \delta$, and the other $\varepsilon^3 \oplus \delta$.
This strongly suggests that these simultaneous eigenspaces are
different Eisenstein lifts of~$z$ from parabolic subgroups of type
$(3,1)$.

\subsection{} 

The Galois representation we find could be an impostor.  There could
be, for example, some irreducible four-di\-men\-sion\-al~$\rho$ that
gives conjugate matrices to ours when evaluated at $\Frob_{\ell}$ for
the few~$\ell$ we can compute and which is the truly attached one.
However, this seems very unlikely.

For~$N$ larger than~$41$, there will be truly attached $\rho$ that our
Galois Finder has not been designed to find.  This happened for
trivial nebentype character in our previous papers \cite{AGM2,AGM3},
where we found lifts of forms from $\mathrm{GSp}(4)$.  This does not
happen in this paper because we are not taking~$N$ big enough.  We
would be delighted to find an apparently attached Galois
representation that is irreducible and not essentially self-dual.
Such a representation would not be a lift from $\mathrm{GSp}(4)$ or
from any proper reductive subgroup of $\SL(4)$.  However, that has not
happened to date.

\subsection{} \label{gal3} We now describe in detail the list of Galois
representations~$\rho$ which our Galois Finder was programmed to use.

We only look at Galois representations whose conductor divides~$N$,
since these are the ones we expect to be constituents of Galois
representations attached to our Hecke eigenclasses of level~$N$.

We begin with Dirichlet characters~$\chi$ with values in~$\F$, which
we identify with one-di\-men\-sion\-al Galois representations as
usual.  We take all the characters of conductor $N_1$ for all $N_1
\mid N$.  Sage's class \texttt{DirichletGroup} enumerates the~$\chi$
automatically.  The characteristic polynomial of Frobenius at~$\ell$
for~$\chi$ is $1 + \chi(\ell) X$, for all $\ell \nmid pN$.
Each~$\chi$ can be lifted to characteristic zero, since $p\equiv 1
\pmod{N}$.

Another one-di\-men\-sion\-al character is the cyclotomic
character~$\varepsilon$.  We look at $\varepsilon^w$ for $w=0,1,2,3$,
because this is predicted by the generalizations of Serre's conjecture
for mod~$p$ Galois representations \cite{AS,ADP}.  These $w$ would
also be the Hodge numbers of the motives conjecturally attached to our
homology eigenclasses.  Our standard list $\mathcal{L}_1$ of
one-di\-men\-sion\-al characters is $\chi \otimes \varepsilon^w$, for
all the~$\chi$ just described and for all $w=0,1,2,3$.

After the Dirichlet characters, we put into the list the Galois
representations~$\rho$ coming from newforms for certain congruence
subgroups of $\SL(2,\Z)$.  We emphasize that these are classical cusp
forms in characteristic zero, even though the~$\rho$ take values in
characteristic $p$.  The characteristic polynomials of Frobenius for
the cusp forms are naturally defined over number fields, so, as we
describe which cusp forms we use, we must also describe how we reduce
to get Galois representations defined over~$\F$.

Let $N_1 \mid N$.  Let $\Q(\zeta_{N_1})$ be the field of $N_1$-th
roots of unity.  Let~$\psi$ be any Dirichlet character of
conductor~$N_1$ taking values in $\C^\times$.  The Galois group
$\Gal(\Q(\zeta_{N_1})/\Q)$ acts on the $\psi$'s by
acting on their values; we take only one~$\psi$ from each Galois
orbit, since the others give Galois-conjugate representations.  Let
$f$ be a newform of weight 2, 3, or 4 for $\Gamma_1(N_1)$ with
nebentype character~$\psi$.  The coefficients of the $q$-expansion
of~$f$ generate a number field~$K_f$, with ring of integers
$\mathcal{O}_{K_f}$.  (This field was called $K_i$ in
Section~\ref{gal2}.)  Let $\mathfrak{P}$ be a prime of~$K_f$ over~$p$.
If $\F$ is of high enough degree over $\F_p$, then the finite field
$\mathcal{O}_{K_f} / \mathfrak{P}$ will have an embedding
$\alpha_{\mathfrak{P}}$ into $\F$.  As we have mentioned, the
extension field~$\F$ of $\F_p$ is chosen so that all these embeddings
will exist.  Therefore, the pair $(f, \mathfrak{P})$ gives rise to a
Galois representation~$\rho$ into $\GL(2, \F)$, by reduction
mod~$\mathfrak{P}$ composed with $\alpha_{\mathfrak{P}}$.  For any
$\ell \nmid pN$, the characteristic polynomial of Frobenius is $1 -
\alpha_{\mathfrak{P}}(a_\ell) X + X^2$, where $a_\ell$ is the
$\ell$-th coefficient in the $q$-expansion of~$f$.  If we chose a
different prime~$\mathfrak{P}$, we would get a Galois-conjugate
representation.

We make a list $\mathcal{L}_2^0$ containing the
representation~$\rho$ for $(f,\mathfrak{P})$, for all $N_1 \mid N$ and
all newforms $f$ of weight 2, 3, or 4 for $\Gamma_1(N_1)$ and all
nebentypes~$\psi$.  Sage's class \texttt{CuspForms}, with its method
\texttt{newforms}, makes this automatic.

We take all the $\rho$ in $\mathcal{L}_2^0$, and tensor
them in all possible ways with the one-di\-men\-sion\-al representations
from the list $\mathcal{L}_1$ of Dirichlet characters and cyclotomic
character powers.  This list of tensor products is our final list
$\mathcal{L}_2$ of two-di\-men\-sion\-al Galois representations.

Our list of three-di\-men\-sion\-al Galois representations is the list of
symmetric squares of $\rho\in\mathcal{L}_2^0$, tensored
in all possible ways with $\mathcal{L}_1$. 

In the results, our Galois representations have a term
$\mathrm{Sym}^2(\sigma)$ for a cusp form~$\sigma$ for three levels,
the prime levels $N=29$, $37$, and $41$.  It all three cases, it
occurs only when~$\eta$ is a quadratic character.  This is because the
symmetric square of the Galois representation attached to a cusp form
with quadratic nebentype can have prime level.

We define the \emph{Hodge-Tate (HT) numbers} for~$\rho$ as follows.
For an element $\chi\otimes \varepsilon^w \in \mathcal{L}_1$, there is
a list of one $\varepsilon$ power, $[w]$.  To a representation coming
from a newform~$\rho$ of weight~$k$, there is a list of two
$\varepsilon$ powers, $[0, k-1]$.  For $\chi\otimes
\varepsilon^w\otimes\rho$, the list is $[w, w+k-1]$.  For direct sums
of representations, the lists are concatenated.  For the
four-di\-men\-sion\-al Galois representations we find to fit our data,
we always observe that the list is $[0,1,2,3]$ after sorting.  This is
what we expected, based on the Serre-type conjectures and the
conjectural HT numbers.  This gives us a check on our computations.
See also Section \ref{HT}.

Another check on our computations comes from considering the
relationship between the nebentype character and the determinant of
the apparently attached representation.  For example, consider a
Galois representation $\rho$ apparently attached to one of our Hecke
eigenclasses, where $\rho$ has the form $\varepsilon^a \oplus \chi
\varepsilon^b \oplus \sigma$, and $\sigma$ is attached to a cusp form
of weight $k$ with nebentype character $\psi$.  Then the determinant
of $\rho$ is $\varepsilon^{a+b+k-1}\chi\psi$ and by the definition of
attachment this must equal $\varepsilon^6\eta$.

\subsection{} \label{galmult} As we have indicated, the Galois groups
$\Gal(\Q(\zeta_N)/\Q)$ or $\Gal(\Q(\zeta_{N_1})/\Q)$ act on our lists
$\mathcal{L}_1$ and $\mathcal{L}_2^0$.  Sometimes a
cohomology group will contain Hecke eigenspaces $E_{(1)}, \dots,
E_{(g)}$ which seem to be attached to Galois representations
$\rho_{(1)}, \dots, \rho_{(g)}$ where $\rho_{(1)}, \dots, \rho_{(g)}$
is an orbit under the Galois action.  We define the \emph{Galois
  multiplicity} of each of $E_{(1)}, \dots, E_{(g)}$ to be~$g$ in this
case.  In the tables in Section~\ref{resMain}, we only list one of the
$E_{(i)}$, and we indicate the Galois multiplicity in the first
column.

In the table in Section~\ref{resMain} for level $N=23$ and $\eta=1$,
for example, we read that~$H^5$ is the direct sum of five
one-dim\-en\-sional Hecke eigenspaces (lines).  The first and second
lines are Galois conjugate, the third and fourth are Galois conjugate,
and the fifth is fixed by Galois.  The Galois multiplicities are
therefore 2, 2, and~1.  Lines one through four are for the cusp
form~$\sigma_{23,2a}$, which is defined over a number field with Galois
group~$G$ of order~2.  Because~$G$ fixes the trivial nebentype
$\eta=1$, it acts on the cohomology.  It interchanges the pairs of
lines 1--2 and 3--4.  The fifth line is for the cusp
form~$\sigma_{23,4}$, which is defined over~$\Q$.  Hence Galois acts
trivially on the fifth line.

\section{Observed regularities in the data and heuristics}\label{ChiReasons}

This section details the regularities we observed in the tables below.
When we have a reasonable heuristic explanation of a pattern, we give
it.  Converting any of these regularities or heuristics to theorems
would require a finer analysis of the Borel-Serre boundary than is
presently available and a greater expertise with Eisenstein series
than we possess.

In this section, we let $\Gamma_0(a,b)$ denote the subgroup of
$\GL(a,\Z)$ where the bottom row is congruent to $(0,\dots,0,*)$
modulo $b$.  Thus $\Gamma_0(N)=\Gamma_0(4,N)\cap \SL(4,\Z)$ in our
notation.  We shall refer to a Hecke eigenclass in
$H^5(\Gamma_0(N),\F_\eta)$ by the letter $z$ and to its attached
Galois representation by~$\rho$.

One pattern mentioned in the previous section is that the determinant
of~$\rho$ always equals~$\varepsilon^6 \eta$.  This is a tautology
from the definition of attachment.

Another pattern we observe is that~$\rho$ must be odd. In other words,
the eigenvalues of $\rho(c)$ are $+1,-1,+1,-1$, where~$c$ denotes
complex conjugation.  This must be the case, as follows from a theorem
of Caraiani and LeHung \cite{CLH}.

One question is why the weights $2,3,4$ occur for the
2-di\-men\-sion\-al irreducible components of the Galois
representations.  Heuristically, the observed weights can be explained
in terms of the homology of the Borel-Serre boundary.  This is
outlined in detail in our first paper \cite{AGM1}, to which we refer
the reader.  Another question is why the exponents of the powers
of~$\varepsilon$ that occur as factors of the 1-di\-men\-sion\-al
components are always contained in the set $\{0,1,2,3\}$, and what is
the relationship between these exponents and the other components.
The heuristic for this comes from deep (conjectural) connections
between Hodge-Tate numbers of Galois representations and the
coefficients of the cohomology classes to which they are attached.
Finally, why do we sometimes observe that the multiplicity of a Hecke
eigenspace is~3, whereas usually it is~1?  This has to do with
oldforms versus newforms.  We now explain these answers in more
detail.

Let $\Gamma=\Gamma_0(N)$.  Let $B_\Gamma$ be the Borel-Serre boundary
of the locally symmetric space $X_\Gamma = \Gamma\backslash
\SL(4,\R)/\mathrm{SO}(4)$. Then $B$ is the union of faces $F(P)$,
where $P$ runs over a set of representatives of $\Gamma$-orbits of
parabolic subgroups $P$ of $\GL(4,\Q)$.  It is simpler to discuss
homology rather than cohomology; this changes nothing qualitatively
about the Hecke eigenvalues and attached Galois representations.

The injection $B_\Gamma\to X_\Gamma \cup B_{\Gamma }$ induces a map on
homology, for any coefficient system~$M$:
\[
H_5(B_\Gamma,M)\to H_5(X_\Gamma \cup B_{\Gamma},M)=H_5(\Gamma,M).
\]
The \emph{boundary homology} is the image of this map.  In this paper,
every class we computed appears to be in the boundary homology.

For each parabolic subgroup~$P$, let $P=LU$, where $L$ is a Levi
component of $P$ and $U$ is the unipotent radical of $P$.  Let
$\pi:P\to P/U$ be the projection.  The image of $\pi$ is isomorphic to
$L$ and is a product of ``blocks" $\GL(n_i,\Q)$, where $\sum n_i = 4$.
If $P$ is conjugate to a standard parabolic subgroup (i.e., one
containing the upper triangular matrices), the block sizes down the
diagonal can be recorded as $(n_1,\dots,n_{k+2})$.  We call this tuple
the ``type" of $P$.  The nonnegative integer $k$ equals the
codimension of $F(P)$ in $B_\Gamma$.  (If there is more than one
parabolic subgroup in the associate class of $P$ we choose the type of
one of them to be the type of $P$.  The ambiguity has no importance
for us.) Below, $k=0$ except in Section \ref{w4}, where $k=1$.

Let $X_L$ denote the symmetric space of $L(\R)$.  Let
$P_\Gamma=P\cap\Gamma$, $U_\Gamma=U\cap\Gamma$, and
$L_\Gamma=\pi(P_\Gamma$).  The face $F(P)$ is a fibration with base
$X_L/L_\Gamma$ and fiber $U(\R)/U_\Gamma$.  The Serre spectral
sequence of this fibration degenerates at $E^2$.  Therefore, if we put
a homology class on each block of $L_\Gamma$, whose degrees
$i_1,\dots,i_{k+2}$ add to $i$, with coefficients in
$H_j(U_\Gamma,M)$, we obtain a class in $H_{i+j}(F(P),M)$.  This class
may or may not give rise to a nonzero class $H_{i+j+k}(B_\Gamma,M)$, depending
on how it behaves in the Leray spectral sequence for the covering of
$B_\Gamma$ by its faces.  Finally, if there is a nonzero class in
$H_{i+j+k}(B_\Gamma,M)$ obtained this way, it may or may not map to a
nonzero class in $H_5(\Gamma,M)$.  All this behaves Hecke-equivariantly.

In this way we expect various kinds of homology Hecke eigenclasses in
the boundary homology, with attached Galois representations that are
reducible, with components corresponding to the homology classes on
the blocks of $L$.  We always have $i+j+k=5$ because we computed
$H_5(\Gamma, \F_\eta)$.

\begin{figure}[htb]
\psfrag{1}{$\scriptstyle 1$}
\psfrag{2}{$\scriptstyle 2$}
\psfrag{3}{$\scriptstyle 3$}
\psfrag{4}{$\scriptstyle 4$}
\psfrag{cusp}{$\scriptstyle \text{cusp}$}
\psfrag{eis}{$\scriptstyle \text{eis}$}
\centering
\subfigure[\label{fig:1a}]{\includegraphics[scale=0.20]{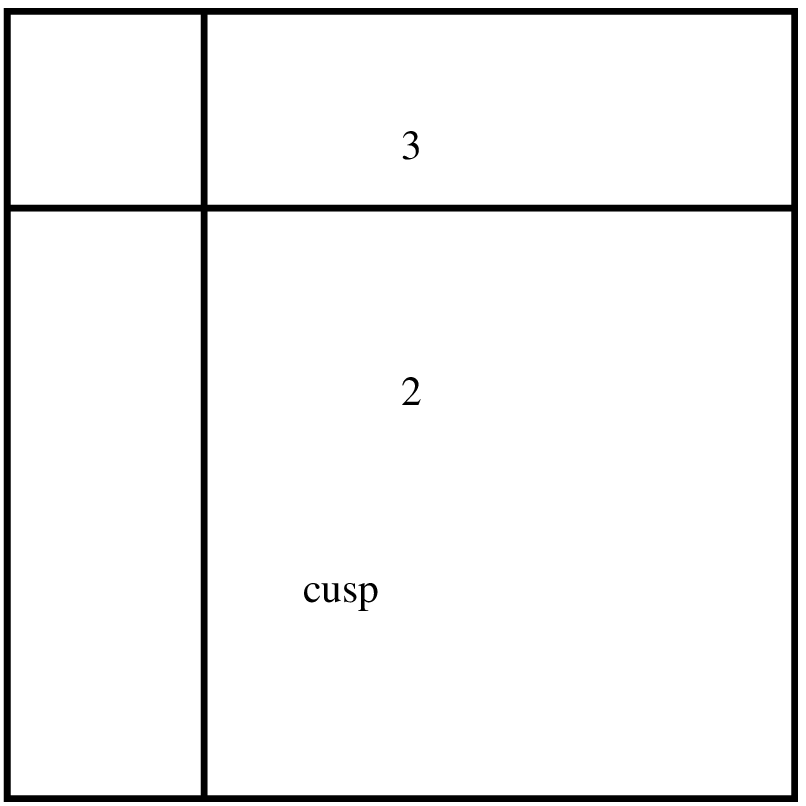}}
\quad\quad
\subfigure[\label{fig:1b}]{\includegraphics[scale=0.20]{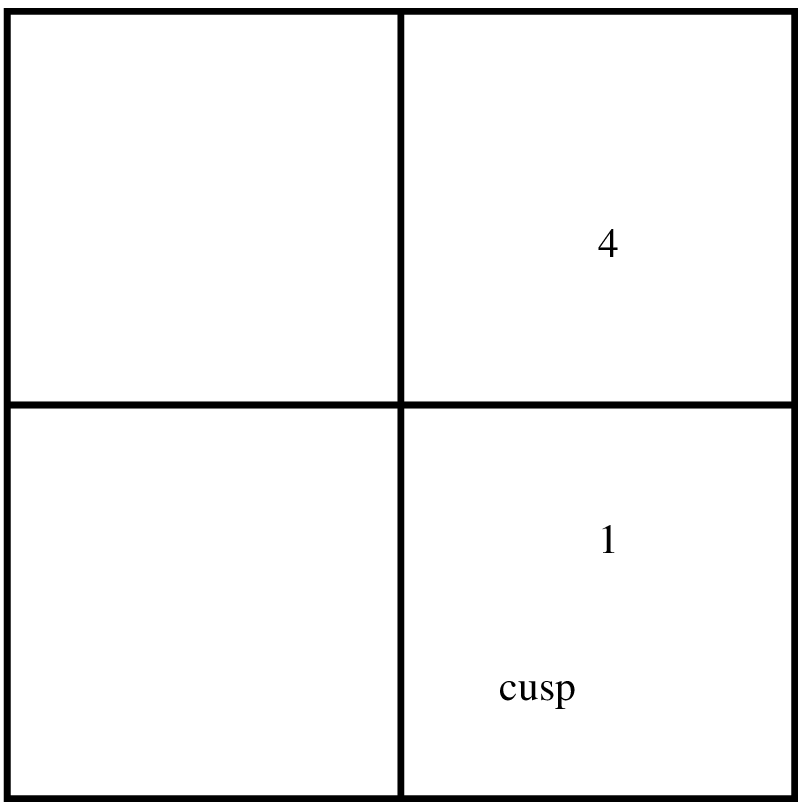}}
\quad\quad 
\subfigure[\label{fig:1c}]{\includegraphics[scale=0.20]{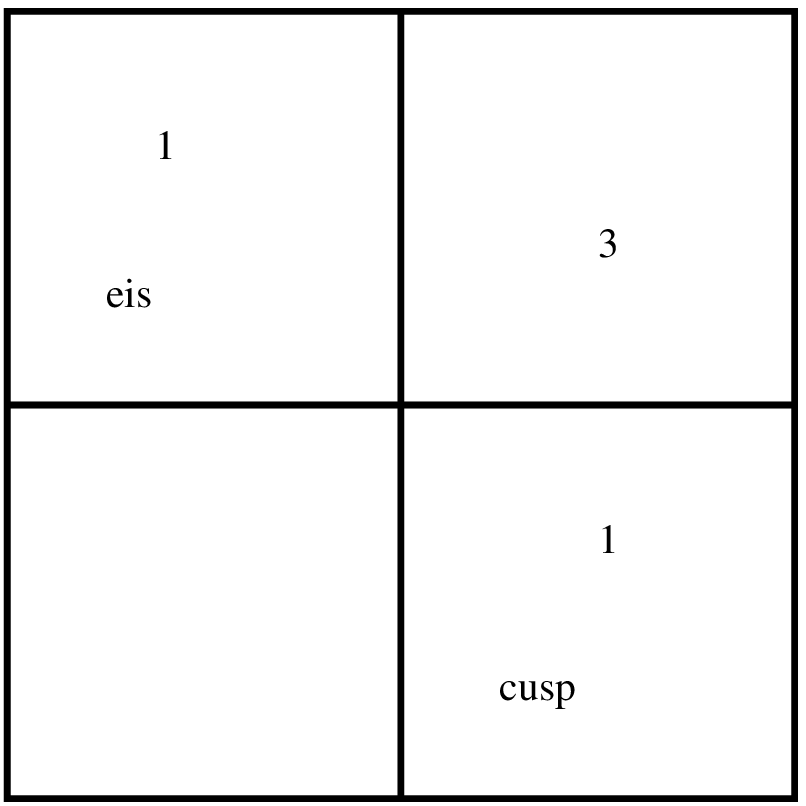}}
\quad\quad
\subfigure[\label{fig:1d}]{\includegraphics[scale=0.20]{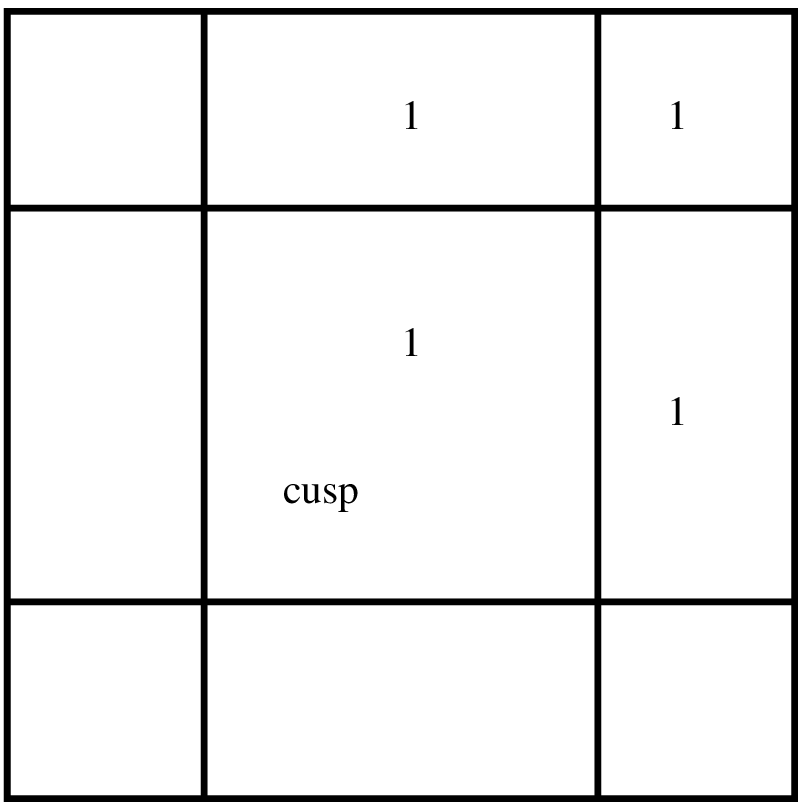}}
\quad\quad 
\subfigure[\label{fig:1e}]{\includegraphics[scale=0.20]{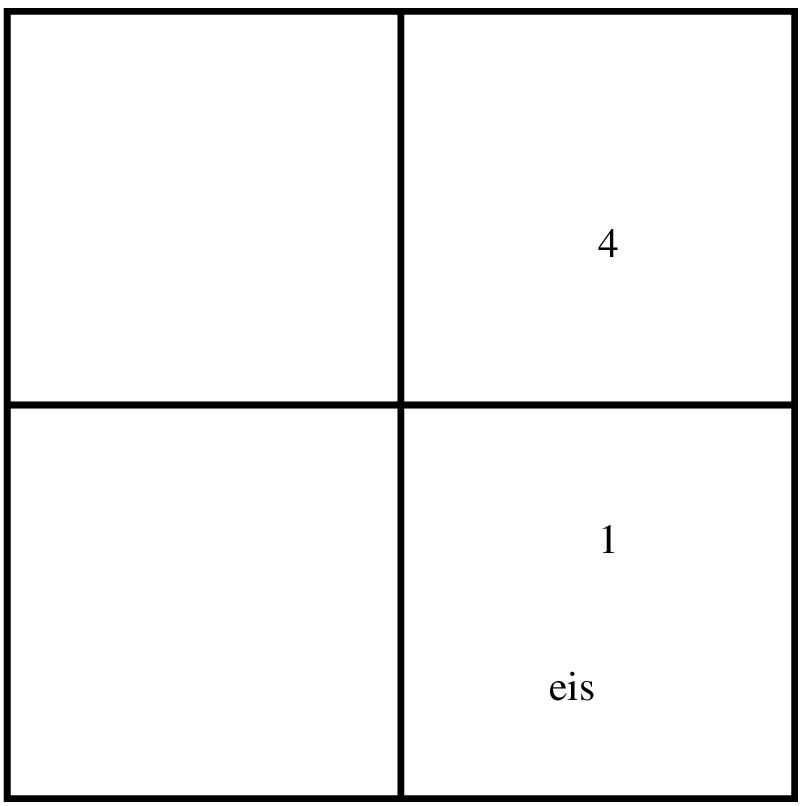}}
\caption{\label{fig:specseq} Schematics of homology classes on faces
  of the Borel-Serre boundary}
\end{figure}

For each type of~$z$ it is convenient to have a schematic picture of
the parameters, as in Figure~\ref{fig:specseq}.  Each
diagram represents a standard parabolic subgroup conjugate to a $P$
that gives rise to some kind of boundary homology.

For a $2\times 2$ block $L'$ of $L$, we use the Eichler-Shimura
theorem to interpret the homology of a congruence subgroup of $L'$
with coefficients in $\mathrm{Sym}^g(\F^2)\otimes \eta$ in terms of
classical modular forms of weight $g+2$ and nebentype~$\eta$.
Therefore the corresponding component of $\rho$ will be attached to
such a modular form.  If we put on $L'$ a class in $H_0(\Gamma\cap
L',M)$, the corresponding component of $\rho$ is observed always to be
the sum of two consecutive powers of the cyclotomic character.  If we
take $H_1(\Gamma\cap L',M)$, we observe either a sum of two characters
(corresponding to an Eisenstein series) or the Galois representation
attached to a cusp form.

A general remark on nebentypes: different $\Gamma$-orbits of the same
type of parabolic subgroups may result in different levels of the
components of $L_\Gamma$.  If $N$ is composite, various nebentype
characters can occur, but they will all have conductor dividing~$N$.

In the following sections we give heuristics along the lines sketched
above that account for all our data.  We reiterate that additional
heuristic schemes are possible and would be needed if we had been
able to push our computations to much higher levels.

\subsection{$\GL(3)$ classes}\label{gl3}

In this case (Figure~1(a)), $P$ is a $(1,3)$-parabolic subgroup;
$i_1=0, i_2=2, j=3$. Note that $H_3(U_\Gamma,\F_\eta)$ is a
one-di\-men\-sion\-al $L'$-module.  We place a cuspidal homology class
$w$ from $H_2(\Gamma_0(3,N),\F_\eta)$ on the second block.  This class
$w$ can be the symmetric square of a classical cusp form, or a class
that is not a symmetric square.  The latter occurs in our data only at
level~$41$.

When $w$ is a symmetric square of the cusp form $s$, the level of $s$
equals $N$, the nebentype of $s$ equals the nebentype $\eta$, and
$\eta$ is the quadratic character.  This is necessary for a symmetric
square at prime level $N$ to have the same level $N$ as the cusp form.

Writing the symmetric square of the Galois representation attached to $w$
as $\tau$, it always appears twice in our data, as $\rho=
\varepsilon^0\oplus \varepsilon\tau$ and $\rho= \varepsilon^3\oplus
\tau$.  This is because there will be two relevant $\Gamma$-orbits of
$P$, corresponding to block sizes $(1,3)$ and $(3,1)$ down the
diagonal.

\subsection{Holomorphic cusp forms of weight 2}\label{w2}

In this case (Figure~1(b)), $P$ is a $(2,2)$-parabolic subgroup;
$i_1=0, i_2=1, j=4$. Note that $H_4(U_\Gamma,\F_\eta)$ is a
one-di\-men\-sion\-al $L'$-module.  We place a cusp form $v$ on one of
the two blocks.  We observe that $v$ always has level~$N$.

In our data, $\sigma$ always appears twice: once in
$\varepsilon^0\oplus \varepsilon^1 \oplus \varepsilon^2 \chi\sigma$
and once in $\varepsilon^2\oplus \varepsilon^3 \oplus \varepsilon^0
\chi\sigma$, for some character $\chi$.  This is because there will be
two relevant $\Gamma$-orbits of $P$, both corresponding to block sizes
$(2,2)$; but in the second orbit, $v$ gets placed on the first block
instead of the second block.

The auxiliary character $\chi$ is the same in both expressions.  We
can and do always choose the ideal $\mathfrak{P}$ so that $\chi=1$.
Sometimes these Galois representations appear with multiplicity~1, and
sometimes with higher multiplicity, and we don't know why the
variability occurs.

\subsection{Holomorphic cusp forms of weight 3}\label{w3}

Cusp forms of odd weight can appear only if $p=2$, as in \cite{AGM6},
or if odd nebentypes are available, as in the current paper.

In this case (Figure~1(c)), $P$ is a $(2,2)$-parabolic subgroup;
$i_1=1, i_2=1, j=3$. Note that $H_3(U_\Gamma,\F_\eta)$ restricted to
either of the $2\times 2$-blocks is a sum of two copies of the
standard 2-di\-men\-sion\-al $\GL(2)$-representation.  We place a
cusp form $v$ of weight 3 on one of the two blocks and an Eisenstein
series $u$ on the other block.

Let $\sigma$ be the Galois representation attached to $v$.  We observe
that $v$ always has level strictly dividing $N$ and always appears in
our data four times as follows:

$\rho= \psi \varepsilon^0\oplus \varepsilon^2 \oplus \varepsilon\sigma$

$\rho=  \varepsilon^0\oplus \psi\varepsilon^2 \oplus \varepsilon\sigma$

$\rho=  \psi\varepsilon^1\oplus \varepsilon^3 \oplus \sigma$ 

$\rho=  \varepsilon^1\oplus \psi\varepsilon^3 \oplus \sigma$

\noindent
with the same character $\psi$ all four times. 

We only have three examples of this, at levels $N=24, 27$ and $28$.
It doesn't always occur even when $N$ is composite and there is an
appropriate $v$ available.  For example, there is a weight 3 cusp form
of level 7 that contributes to $N=28$ but it does not contribute when
$N=14$.  We have no conjecture as to when a weight~3 cusp form appears
for a given $(N,\eta)$.

\subsection{Holomorphic cusp forms of weight 4}\label{w4}

In this case (Figure~1(d)), $P$ is a $(1,2,1)$-parabolic subgroup;
$i_1=0, i_2=1, i_3=0, j=3$. Note that $H_3(U_\Gamma,\F_\eta)$ contains
an $L'$-submodule isomorphic to $\mathrm{Sym}^2$ of the standard
representation.  We place a cusp form $v$ of weight $4$ on the second
block.

Let $\sigma$ be the Galois representation attached to $v$.  We observe
that $\rho= \varepsilon^1\oplus \varepsilon^2 \oplus \sigma$ occurs
only once in our data, if at all.  It occurs if and only if the
special value $L(v,1/2)$ of the $L$-function is~$0$.  For the levels
we have computed, this occurs only when $\eta=1$.  The level of $v$
always divides $N$ but need not equal $N$.

\subsection{Sums of 4 characters}\label{allchars}

See Figure~1(e).  Here, as in~(\ref{w2}), $P$ is a $(2,2)$-parabolic
subgroup; $i_1=0, i_2=1, j=4$. We place an Eisenstein series $e$ on
one of the two blocks.  Not surprisingly, $e$ always has level
dividing $N$ and the two characters $\psi$ and $\chi$ associated with
$e$ have conductors dividing $N$.

The following behavior is mysterious to us.    If $\eta$ factors nontrivially as 
$\eta=\psi\chi$ then either all three of the following or none of the following occur:

$\rho=  \psi\varepsilon^0\oplus \chi \varepsilon^1\oplus\varepsilon^2\oplus \varepsilon^3$

$\rho= \psi \varepsilon^0\oplus \varepsilon^1\oplus\varepsilon^2\oplus \chi\varepsilon^3$

$\rho=  \varepsilon^0\oplus \varepsilon^1\oplus\psi\varepsilon^2\oplus \chi\varepsilon^3$

\noindent For example, when $N=9$ all three forms occur, and when
$N=13$ none of the three occur.  Note that in a given triple of
$\rho$'s, one of the characters multiplies even powers of
$\varepsilon$ and the other multiplies odd powers, which to some
extent is explained because $\rho$ must be odd.  But why don't we ever
get $\rho= \varepsilon^0\oplus
\chi\varepsilon^1\oplus\psi\varepsilon^2\oplus \varepsilon^3$?  For
some reason, the unadorned powers of $\varepsilon$ are always
consecutive.  This is true of all the patterns observed above, except
for those in~(\ref{w3}).  In the case of weight 3, the two
1-di\-men\-sion\-al components of $\rho$ are not both naked powers of
$\varepsilon$.

Note that $\psi$ and $\chi$ can trade places to get another triple,
giving 6 $\rho$'s in total, for example, when $n=15$.  Factoring of
$\eta$ seems to be important here.  For example, when $\eta=1$ we
never get $\rho= \varepsilon^0\oplus
\varepsilon^1\oplus\varepsilon^2\oplus \varepsilon^3$.  However, for
$\GL(n)$ with $n$ larger than $4$, the sum of consecutive powers of
$\varepsilon$ may be attached to a homology class, for example in the
case of a Borel stable class \cite[(4.2)]{ashconjec}.

\subsection{Powers of $\varepsilon$}\label{HT}
Another observed pattern has to do with the powers of the cyclotomic
character that appear in $\rho$.  Let us say that $\varepsilon^i$ has
HT (Hodge-Tate) number $i$.  We assign HT numbers $0,k-1$ to a Galois
representation attached to a cusp form of weight $k$.  We assign HT
numbers $0,k-1,2(k-1)$ to a Galois representation attached to the
symmetric square of a cusp form of weight $k$.  (In our data, we only
see symmetric squares when $k=2$.)  We assign HT numbers $0,1,2$ to a
Galois representation attached to a $\GL(3)$-homology eigenclass with
trivial coefficients.  If a Galois representation is tensored with
$\varepsilon^i$, then add $i$ to each of its HT-numbers.

A folklore conjecture in the theory of arithmetic cohomology predicts
that a Galois representation attached to a Hecke eigenclass in
$H^5(\Gamma, \F_\eta)$ should have HT numbers $0,1,2,3$.  This is
observed in all of our data.

\subsection{Hecke multiplicity 3}\label{mult=3}

We defined \emph{Hecke multiplicity} in Section~\ref{gal1}.  In every
case of our data, the Hecke multiplicity of the eigenspace for a
system of Hecke eigenvalues equals either 1 or 3.  As stated in
Section \ref{gal}, the Hecke operators we computed are always observed
to be semisimple.

Hecke multiplicity 3 occurs in our data only when $N$ is composite and
the components of $\rho$ have conductors strictly dividing $N$.  In
our data, this happens for $N=18, 22, 26, 27$, and $28$.  When $N=18$
or $27$ the relevant $\rho$ is a sum of four characters, involving in
two of the summands the quadratic character of conductor~3.  In the
other three cases, one of the components of $\rho$ is attached to a
cusp form of weight 2 and level $N/2$.  We see no general rule as to
why these cases of Hecke multiplicity 3 occur and not others that
would be possible.

Here is a partial explanation of why the multiplicity is 3 rather than
some other number when one of the components of $\rho$ is attached to
a cusp form $\sigma$ of weight 2 and level $N/2$.  First consider a
parabolic subgroup of type $(2,2)$ such that $L_\Gamma$ is isomorphic
to a subgroup of index two in $\Gamma_0(2,N)\times \GL(2,\Z)$. Since
$\sigma$ is an oldform for $\Gamma_0(2,N)$, its system of Hecke
eigenvalues contributes twice to the cohomology of $L_\Gamma$, giving
a Hecke multiplicity of 2, so far, for the system of Hecke eigenvalues
to which $\rho$ is attached.  However, there is another $\Gamma$-orbit
of parabolic subgroups of type $(2,2)$ such that $L_\Gamma$ is
isomorphic to a subgroup of index two in $\Gamma_0(2,N/2)\times
\Gamma_0(2,2)$.  Here, $\sigma$ is a newform for $\Gamma_0(2,N/2)$, so
its system of Hecke eigenvalues contributes once to the cohomology of
$L_\Gamma$, adding~1 to the Hecke multiplicity for the system of Hecke
eigenvalues to which $\rho$ is attached.  The total is $2+1=3$.

\section{Results}\label{res}

\subsection{} \label{resMain}
The tables in this section present the main results of the paper.

The topmost box in each table gives the level~$N$, the nebentype~$\eta$,
and the field~$\F_{p^r} = GF(p^r)$ that was our proxy for~$\C$.  We
only include one representative for each Galois orbit of nebentype
characters.  Next we list the Hecke operators we computed.  $T_\ell$
means we computed $T_{\ell,1}$, $T_{\ell,2}$, and $T_{\ell,3}$.
Listing $T_{\ell,1}$ means we computed only that part of $T_\ell$.

The succeeding rows in each table give the Galois multiplicity
(Section~\ref{galmult}), the Hecke multiplicity (Section~\ref{gal1}),
and the Galois representation itself.

The characters~$\chi_{N}$ or~$\chi_{N,i}$ are a basis for the
Dirichlet characters $(\Z/N\Z)^\times \to \F_p$.  They are listed in a
separate table in Section~\ref{resChi}.  The cyclotomic character is
denoted~$\varepsilon$.

The~$\sigma_{N,k}$ are classical cuspidal homomorphic newforms of
level~$N$ and weight~$k$.  They are listed in a separate table in
Section~\ref{resCusp}.  We use the same symbol $\sigma_{N,k}$ to stand
for the two-di\-men\-sional Galois representation attached to the cusp
form of that name.  When we have more than one cusp form for the
same~$N$ and~$k$, we give them names like $\sigma_{17,2a}$ and
$\sigma_{17,2b}$.  The symmetric square of $\sigma$ is denoted
$\mathrm{Sym}^2(\sigma)$.

The $\SL_3$ representation~$\delta$ is defined in Section~\ref{dihedralAndDoud}.

\begin{center}
\begin{tabular}{|r|r|c|}
\hline
\multicolumn{3}{|l|}{\textbf{Level $N = 9$.  Nebentype $\eta = 1$.  Field $\mathbb{F} = GF(12379)$.}} \\
\multicolumn{3}{|l|}{Computed $T_{2}$, $T_{5}$, $T_{7}$.} \\
\hline
1 & $1$ & $\varepsilon^{0} \oplus \varepsilon^{1} \oplus \chi_{9}^{3} \varepsilon^{2} \oplus \chi_{9}^{3} \varepsilon^{3}$ \\
\hline
1 & $1$ & $\chi_{9}^{3} \varepsilon^{0} \oplus \varepsilon^{1} \oplus \varepsilon^{2} \oplus \chi_{9}^{3} \varepsilon^{3}$ \\
\hline
1 & $1$ & $\chi_{9}^{3} \varepsilon^{0} \oplus \chi_{9}^{3} \varepsilon^{1} \oplus \varepsilon^{2} \oplus \varepsilon^{3}$ \\
\hline
\end{tabular}
\end{center}

\begin{center}
\begin{tabular}{|r|r|c|}
\hline
\multicolumn{3}{|l|}{\textbf{Level $N = 11$.  Nebentype $\eta = 1$.  Field $\mathbb{F} = GF(4001^{2})$.}} \\
\multicolumn{3}{|l|}{Computed $T_{2}$, $T_{3}$, $T_{5}$, $T_{7}$.} \\
\hline
1 & $1$ & $\varepsilon^{0} \oplus \varepsilon^{1} \oplus \varepsilon^{2} \sigma_{11,2}$ \\
\hline
1 & $1$ & $\varepsilon^{2} \oplus \varepsilon^{3} \oplus \varepsilon^{0} \sigma_{11,2}$ \\
\hline
\end{tabular}
\end{center}

\begin{center}
\begin{tabular}{|r|r|c|}
\hline
\multicolumn{3}{|l|}{\textbf{Level $N = 12$.  Nebentype $\eta = \chi_{12,0}\chi_{12,1}$.  Field $\mathbb{F} = GF(5413^{2})$.}} \\
\multicolumn{3}{|l|}{Computed $T_{5}$, $T_{7}$.} \\
\hline
1 & $1$ & $\varepsilon^{0} \oplus \varepsilon^{1} \oplus \chi_{12,0} \varepsilon^{2} \oplus \chi_{12,1} \varepsilon^{3}$ \\
\hline
1 & $1$ & $\varepsilon^{0} \oplus \varepsilon^{1} \oplus \chi_{12,1} \varepsilon^{2} \oplus \chi_{12,0} \varepsilon^{3}$ \\
\hline
1 & $1$ & $\chi_{12,0} \varepsilon^{0} \oplus \varepsilon^{1} \oplus \varepsilon^{2} \oplus \chi_{12,1} \varepsilon^{3}$ \\
\hline
1 & $1$ & $\chi_{12,1} \varepsilon^{0} \oplus \varepsilon^{1} \oplus \varepsilon^{2} \oplus \chi_{12,0} \varepsilon^{3}$ \\
\hline
1 & $1$ & $\chi_{12,0} \varepsilon^{0} \oplus \chi_{12,1} \varepsilon^{1} \oplus \varepsilon^{2} \oplus \varepsilon^{3}$ \\
\hline
1 & $1$ & $\chi_{12,1} \varepsilon^{0} \oplus \chi_{12,0} \varepsilon^{1} \oplus \varepsilon^{2} \oplus \varepsilon^{3}$ \\
\hline
\end{tabular}
\end{center}

\begin{center}
\begin{tabular}{|r|r|c|}
\hline
\multicolumn{3}{|l|}{\textbf{Level $N = 13$.  Nebentype $\eta = 1$.  Field $\mathbb{F} = GF(12037)$.}} \\
\multicolumn{3}{|l|}{Computed $T_{2}$, $T_{3}$, $T_{5}$, $T_{7}$.} \\
\hline
1 & $1$ & $\varepsilon^{1} \oplus \varepsilon^{2} \oplus \varepsilon^{0} \sigma_{13,4}$ \\
\hline
\end{tabular}
\end{center}

\begin{center}
\begin{tabular}{|r|r|c|}
\hline
\multicolumn{3}{|l|}{\textbf{Level $N = 13$.  Nebentype $\eta = \chi_{13}^{2}$.  Field $\mathbb{F} = GF(12037)$.}} \\
\multicolumn{3}{|l|}{Computed $T_{2}$, $T_{3}$, $T_{5}$, $T_{7}$.} \\
\hline
1 & $1$ & $\varepsilon^{0} \oplus \varepsilon^{1} \oplus \varepsilon^{2} \sigma_{13,2}$ \\
\hline
1 & $1$ & $\varepsilon^{2} \oplus \varepsilon^{3} \oplus \varepsilon^{0} \sigma_{13,2}$ \\
\hline
\end{tabular}
\end{center}

\begin{center}
\begin{tabular}{|r|r|c|}
\hline
\multicolumn{3}{|l|}{\textbf{Level $N = 14$.  Nebentype $\eta = 1$.  Field $\mathbb{F} = GF(12379^{2})$.}} \\
\multicolumn{3}{|l|}{Computed $T_{3}$, $T_{5}$.} \\
\hline
1 & $1$ & $\varepsilon^{0} \oplus \varepsilon^{1} \oplus \varepsilon^{2} \sigma_{14,2}$ \\
\hline
1 & $1$ & $\varepsilon^{2} \oplus \varepsilon^{3} \oplus \varepsilon^{0} \sigma_{14,2}$ \\
\hline
\end{tabular}
\end{center}

\begin{center}
\begin{tabular}{|r|r|c|}
\hline
\multicolumn{3}{|l|}{\textbf{Level $N = 15$.  Nebentype $\eta = 1$.  Field $\mathbb{F} = GF(12037^{2})$.}} \\
\multicolumn{3}{|l|}{Computed $T_{2}$, $T_{7}$, $T_{11,1}$.} \\
\hline
1 & $1$ & $\varepsilon^{0} \oplus \varepsilon^{1} \oplus \varepsilon^{2} \sigma_{15,2}$ \\
\hline
1 & $1$ & $\varepsilon^{2} \oplus \varepsilon^{3} \oplus \varepsilon^{0} \sigma_{15,2}$ \\
\hline
\end{tabular}
\end{center}

\begin{center}
\begin{tabular}{|r|r|c|}
\hline
\multicolumn{3}{|l|}{\textbf{Level $N = 15$.  Nebentype $\eta = \chi_{15,0}\chi_{15,1}$.  Field $\mathbb{F} = GF(12037^{2})$.}} \\
\multicolumn{3}{|l|}{Computed $T_{2}$, $T_{7}$.} \\
\hline
1 & $1$ & $\varepsilon^{0} \oplus \varepsilon^{1} \oplus \chi_{15,0} \varepsilon^{2} \oplus \chi_{15,1} \varepsilon^{3}$ \\
\hline
1 & $1$ & $\varepsilon^{0} \oplus \varepsilon^{1} \oplus \chi_{15,1} \varepsilon^{2} \oplus \chi_{15,0} \varepsilon^{3}$ \\
\hline
1 & $1$ & $\chi_{15,0} \varepsilon^{0} \oplus \varepsilon^{1} \oplus \varepsilon^{2} \oplus \chi_{15,1} \varepsilon^{3}$ \\
\hline
1 & $1$ & $\chi_{15,1} \varepsilon^{0} \oplus \varepsilon^{1} \oplus \varepsilon^{2} \oplus \chi_{15,0} \varepsilon^{3}$ \\
\hline
1 & $1$ & $\chi_{15,0} \varepsilon^{0} \oplus \chi_{15,1} \varepsilon^{1} \oplus \varepsilon^{2} \oplus \varepsilon^{3}$ \\
\hline
1 & $1$ & $\chi_{15,1} \varepsilon^{0} \oplus \chi_{15,0} \varepsilon^{1} \oplus \varepsilon^{2} \oplus \varepsilon^{3}$ \\
\hline
\end{tabular}
\end{center}

\begin{center}
\begin{tabular}{|r|r|c|}
\hline
\multicolumn{3}{|l|}{\textbf{Level $N = 16$.  Nebentype $\eta = 1$.  Field $\mathbb{F} = GF(4001^{6})$.}} \\
\multicolumn{3}{|l|}{Computed $T_{3}$, $T_{5}$, $T_{7}$.} \\
\hline
1 & $1$ & $\varepsilon^{0} \oplus \varepsilon^{1} \oplus \chi_{16,0} \varepsilon^{2} \oplus \chi_{16,0} \varepsilon^{3}$ \\
\hline
1 & $1$ & $\chi_{16,0} \varepsilon^{0} \oplus \varepsilon^{1} \oplus \varepsilon^{2} \oplus \chi_{16,0} \varepsilon^{3}$ \\
\hline
1 & $1$ & $\chi_{16,0} \varepsilon^{0} \oplus \chi_{16,0} \varepsilon^{1} \oplus \varepsilon^{2} \oplus \varepsilon^{3}$ \\
\hline
\end{tabular}
\end{center}

\begin{center}
\begin{tabular}{|r|r|c|}
\hline
\multicolumn{3}{|l|}{\textbf{Level $N = 16$.  Nebentype $\eta = \chi_{16,1}$.  Field $\mathbb{F} = GF(4001^{6})$.}} \\
\multicolumn{3}{|l|}{Computed $T_{3}$, $T_{5}$, $T_{7}$.} \\
\hline
1 & $1$ & $\varepsilon^{0} \oplus \varepsilon^{1} \oplus \varepsilon^{2} \sigma_{16,2}$ \\
\hline
1 & $1$ & $\varepsilon^{2} \oplus \varepsilon^{3} \oplus \varepsilon^{0} \sigma_{16,2}$ \\
\hline
\end{tabular}
\end{center}

\begin{center}
\begin{tabular}{|r|r|c|}
\hline
\multicolumn{3}{|l|}{\textbf{Level $N = 17$.  Nebentype $\eta = 1$.  Field $\mathbb{F} = GF(16001^{2})$.}} \\
\multicolumn{3}{|l|}{Computed $T_{2}$, $T_{3}$, $T_{5}$, $T_{7}$.} \\
\hline
1 & $1$ & $\varepsilon^{0} \oplus \varepsilon^{1} \oplus \varepsilon^{2} \sigma_{17,2a}$ \\
\hline
1 & $1$ & $\varepsilon^{2} \oplus \varepsilon^{3} \oplus \varepsilon^{0} \sigma_{17,2a}$ \\
\hline
1 & $1$ & $\varepsilon^{1} \oplus \varepsilon^{2} \oplus \varepsilon^{0} \sigma_{17,4}$ \\
\hline
\end{tabular}
\end{center}

\begin{center}
\begin{tabular}{|r|r|c|}
\hline
\multicolumn{3}{|l|}{\textbf{Level $N = 17$.  Nebentype $\eta = \chi_{17}^{2}$.  Field $\mathbb{F} = GF(16001^{2})$.}} \\
\multicolumn{3}{|l|}{Computed $T_{2}$, $T_{3}$, $T_{5}$, $T_{7}$.} \\
\hline
1 & $1$ & $\varepsilon^{0} \oplus \varepsilon^{1} \oplus \varepsilon^{2} \sigma_{17,2b}$ \\
\hline
1 & $1$ & $\varepsilon^{2} \oplus \varepsilon^{3} \oplus \varepsilon^{0} \sigma_{17,2b}$ \\
\hline
\end{tabular}
\end{center}

\begin{center}
\begin{tabular}{|r|r|c|}
\hline
\multicolumn{3}{|l|}{\textbf{Level $N = 18$.  Nebentype $\eta = 1$.  Field $\mathbb{F} = GF(3637^{2})$.}} \\
\multicolumn{3}{|l|}{Computed $T_{5}$, $T_{7}$.} \\
\hline
1 & $3$ & $\varepsilon^{0} \oplus \varepsilon^{1} \oplus \chi_{18}^{3} \varepsilon^{2} \oplus \chi_{18}^{3} \varepsilon^{3}$ \\
\hline
1 & $3$ & $\chi_{18}^{3} \varepsilon^{0} \oplus \varepsilon^{1} \oplus \varepsilon^{2} \oplus \chi_{18}^{3} \varepsilon^{3}$ \\
\hline
1 & $3$ & $\chi_{18}^{3} \varepsilon^{0} \oplus \chi_{18}^{3} \varepsilon^{1} \oplus \varepsilon^{2} \oplus \varepsilon^{3}$ \\
\hline
\end{tabular}
\end{center}

\begin{center}
\begin{tabular}{|r|r|c|}
\hline
\multicolumn{3}{|l|}{\textbf{Level $N = 18$.  Nebentype $\eta = \chi_{18}^{2}$.  Field $\mathbb{F} = GF(3637^{2})$.}} \\
\multicolumn{3}{|l|}{Computed $T_{5}$, $T_{7}$, $T_{11,1}$.} \\
\hline
1 & $1$ & $\varepsilon^{0} \oplus \varepsilon^{1} \oplus \varepsilon^{2} \sigma_{18,2}$ \\
\hline
1 & $1$ & $\varepsilon^{2} \oplus \varepsilon^{3} \oplus \varepsilon^{0} \sigma_{18,2}$ \\
\hline
\end{tabular}
\end{center}

\begin{center}
\begin{tabular}{|r|r|c|}
\hline
\multicolumn{3}{|l|}{\textbf{Level $N = 19$.  Nebentype $\eta = 1$.  Field $\mathbb{F} = GF(3637^{6})$.}} \\
\multicolumn{3}{|l|}{Computed $T_{2}$, $T_{3}$, $T_{5}$, $T_{7}$.} \\
\hline
1 & $1$ & $\varepsilon^{0} \oplus \varepsilon^{1} \oplus \varepsilon^{2} \sigma_{19,2a}$ \\
\hline
1 & $1$ & $\varepsilon^{2} \oplus \varepsilon^{3} \oplus \varepsilon^{0} \sigma_{19,2a}$ \\
\hline
1 & $1$ & $\varepsilon^{1} \oplus \varepsilon^{2} \oplus \varepsilon^{0} \sigma_{19,4}$ \\
\hline
\end{tabular}
\end{center}

\begin{center}
\begin{tabular}{|r|r|c|}
\hline
\multicolumn{3}{|l|}{\textbf{Level $N = 19$.  Nebentype $\eta = \chi_{19}^{2}$.  Field $\mathbb{F} = GF(3637^{6})$.}} \\
\multicolumn{3}{|l|}{Computed $T_{2}$, $T_{3}$, $T_{5}$, $T_{7}$.} \\
\hline
1 & $1$ & $\varepsilon^{0} \oplus \varepsilon^{1} \oplus \varepsilon^{2} \sigma_{19,2b}$ \\
\hline
1 & $1$ & $\varepsilon^{2} \oplus \varepsilon^{3} \oplus \varepsilon^{0} \sigma_{19,2b}$ \\
\hline
\end{tabular}
\end{center}

\begin{center}
\begin{tabular}{|r|r|c|}
\hline
\multicolumn{3}{|l|}{\textbf{Level $N = 20$.  Nebentype $\eta = 1$.  Field $\mathbb{F} = GF(12037^{12})$.}} \\
\multicolumn{3}{|l|}{Computed $T_{3}$, $T_{7}$, $T_{11,1}$, $T_{13,1}$.} \\
\hline
1 & $1$ & $\varepsilon^{0} \oplus \varepsilon^{1} \oplus \varepsilon^{2} \sigma_{20,2a}$ \\
\hline
1 & $1$ & $\varepsilon^{2} \oplus \varepsilon^{3} \oplus \varepsilon^{0} \sigma_{20,2a}$ \\
\hline
\end{tabular}
\end{center}

\begin{center}
\begin{tabular}{|r|r|c|}
\hline
\multicolumn{3}{|l|}{\textbf{Level $N = 20$.  Nebentype $\eta = \chi_{20,0}\chi_{20,1}$.  Field $\mathbb{F} = GF(12037^{12})$.}} \\
\multicolumn{3}{|l|}{Computed $T_{3}$, $T_{7}$, $T_{11,1}$, $T_{13,1}$.} \\
\hline
1 & $1$ & $\varepsilon^{0} \oplus \varepsilon^{1} \oplus \chi_{20,0} \varepsilon^{2} \oplus \chi_{20,1} \varepsilon^{3}$ \\
\hline
1 & $1$ & $\varepsilon^{0} \oplus \varepsilon^{1} \oplus \chi_{20,1} \varepsilon^{2} \oplus \chi_{20,0} \varepsilon^{3}$ \\
\hline
1 & $1$ & $\chi_{20,0} \varepsilon^{0} \oplus \varepsilon^{1} \oplus \varepsilon^{2} \oplus \chi_{20,1} \varepsilon^{3}$ \\
\hline
1 & $1$ & $\chi_{20,1} \varepsilon^{0} \oplus \varepsilon^{1} \oplus \varepsilon^{2} \oplus \chi_{20,0} \varepsilon^{3}$ \\
\hline
1 & $1$ & $\chi_{20,0} \varepsilon^{0} \oplus \chi_{20,1} \varepsilon^{1} \oplus \varepsilon^{2} \oplus \varepsilon^{3}$ \\
\hline
1 & $1$ & $\chi_{20,1} \varepsilon^{0} \oplus \chi_{20,0} \varepsilon^{1} \oplus \varepsilon^{2} \oplus \varepsilon^{3}$ \\
\hline
1 & $1$ & $\varepsilon^{0} \oplus \varepsilon^{1} \oplus \varepsilon^{2} \sigma_{20,2b}$ \\
\hline
1 & $1$ & $\varepsilon^{2} \oplus \varepsilon^{3} \oplus \varepsilon^{0} \sigma_{20,2b}$ \\
\hline
\end{tabular}
\end{center}

\begin{center}
\begin{tabular}{|r|r|c|}
\hline
\multicolumn{3}{|l|}{\textbf{Level $N = 21$.  Nebentype $\eta = 1$.  Field $\mathbb{F} = GF(12037^{6})$.}} \\
\multicolumn{3}{|l|}{Computed $T_{2}$, $T_{5}$.} \\
\hline
1 & $1$ & $\varepsilon^{0} \oplus \varepsilon^{1} \oplus \varepsilon^{2} \sigma_{21,2a}$ \\
\hline
1 & $1$ & $\varepsilon^{2} \oplus \varepsilon^{3} \oplus \varepsilon^{0} \sigma_{21,2a}$ \\
\hline
1 & $1$ & $\varepsilon^{1} \oplus \varepsilon^{2} \oplus \varepsilon^{0} \sigma_{21,4}$ \\
\hline
\end{tabular}
\end{center}

\begin{center}
\begin{tabular}{|r|r|c|}
\hline
\multicolumn{3}{|l|}{\textbf{Level $N = 21$.  Nebentype $\eta = \chi_{21,1}^{2}$.  Field $\mathbb{F} = GF(12037^{6})$.}} \\
\multicolumn{3}{|l|}{Computed $T_{2}$, $T_{5}$.} \\
\hline
1 & $1$ & $\varepsilon^{0} \oplus \varepsilon^{1} \oplus \varepsilon^{2} \sigma_{21,2b}$ \\
\hline
1 & $1$ & $\varepsilon^{2} \oplus \varepsilon^{3} \oplus \varepsilon^{0} \sigma_{21,2b}$ \\
\hline
\end{tabular}
\end{center}

\begin{center}
\begin{tabular}{|r|r|c|}
\hline
\multicolumn{3}{|l|}{\textbf{Level $N = 21$.  Nebentype $\eta = \chi_{21,0}\chi_{21,1}$.  Field $\mathbb{F} = GF(12037^{6})$.}} \\
\multicolumn{3}{|l|}{Computed $T_{2}$, $T_{5}$, $T_{11,1}$, $T_{13,1}$.} \\
\hline
1 & $1$ & $\varepsilon^{0} \oplus \varepsilon^{1} \oplus \chi_{21,0} \varepsilon^{2} \oplus \chi_{21,1} \varepsilon^{3}$ \\
\hline
1 & $1$ & $\varepsilon^{0} \oplus \varepsilon^{1} \oplus \chi_{21,1} \varepsilon^{2} \oplus \chi_{21,0} \varepsilon^{3}$ \\
\hline
1 & $1$ & $\chi_{21,0} \varepsilon^{0} \oplus \varepsilon^{1} \oplus \varepsilon^{2} \oplus \chi_{21,1} \varepsilon^{3}$ \\
\hline
1 & $1$ & $\chi_{21,1} \varepsilon^{0} \oplus \varepsilon^{1} \oplus \varepsilon^{2} \oplus \chi_{21,0} \varepsilon^{3}$ \\
\hline
1 & $1$ & $\chi_{21,0} \varepsilon^{0} \oplus \chi_{21,1} \varepsilon^{1} \oplus \varepsilon^{2} \oplus \varepsilon^{3}$ \\
\hline
1 & $1$ & $\chi_{21,1} \varepsilon^{0} \oplus \chi_{21,0} \varepsilon^{1} \oplus \varepsilon^{2} \oplus \varepsilon^{3}$ \\
\hline
1 & $1$ & $\varepsilon^{0} \oplus \varepsilon^{1} \oplus \varepsilon^{2} \sigma_{21,2c}$ \\
\hline
1 & $1$ & $\varepsilon^{2} \oplus \varepsilon^{3} \oplus \varepsilon^{0} \sigma_{21,2c}$ \\
\hline
\end{tabular}
\end{center}

\begin{center}
\begin{tabular}{|r|r|c|}
\hline
\multicolumn{3}{|l|}{\textbf{Level $N = 21$.  Nebentype $\eta = \chi_{21,0}\chi_{21,1}^{3}$.  Field $\mathbb{F} = GF(12037^{6})$.}} \\
\multicolumn{3}{|l|}{Computed $T_{2}$, $T_{5}$.} \\
\hline
1 & $1$ & $\varepsilon^{0} \oplus \varepsilon^{1} \oplus \chi_{21,0} \varepsilon^{2} \oplus \chi_{21,1}^{3} \varepsilon^{3}$ \\
\hline
1 & $1$ & $\varepsilon^{0} \oplus \varepsilon^{1} \oplus \chi_{21,1}^{3} \varepsilon^{2} \oplus \chi_{21,0} \varepsilon^{3}$ \\
\hline
1 & $1$ & $\chi_{21,0} \varepsilon^{0} \oplus \varepsilon^{1} \oplus \varepsilon^{2} \oplus \chi_{21,1}^{3} \varepsilon^{3}$ \\
\hline
1 & $1$ & $\chi_{21,1}^{3} \varepsilon^{0} \oplus \varepsilon^{1} \oplus \varepsilon^{2} \oplus \chi_{21,0} \varepsilon^{3}$ \\
\hline
1 & $1$ & $\chi_{21,0} \varepsilon^{0} \oplus \chi_{21,1}^{3} \varepsilon^{1} \oplus \varepsilon^{2} \oplus \varepsilon^{3}$ \\
\hline
1 & $1$ & $\chi_{21,1}^{3} \varepsilon^{0} \oplus \chi_{21,0} \varepsilon^{1} \oplus \varepsilon^{2} \oplus \varepsilon^{3}$ \\
\hline
1 & $1$ & $\varepsilon^{0} \oplus \chi_{21,0} \varepsilon^{2} \oplus \varepsilon^{1} \sigma_{7,3}$ \\
\hline
1 & $1$ & $\varepsilon^{1} \oplus \chi_{21,0} \varepsilon^{3} \oplus \varepsilon^{0} \sigma_{7,3}$ \\
\hline
1 & $1$ & $\chi_{21,0} \varepsilon^{0} \oplus \varepsilon^{2} \oplus \varepsilon^{1} \sigma_{7,3}$ \\
\hline
1 & $1$ & $\chi_{21,0} \varepsilon^{1} \oplus \varepsilon^{3} \oplus \varepsilon^{0} \sigma_{7,3}$ \\
\hline
\end{tabular}
\end{center}

\begin{center}
\begin{tabular}{|r|r|c|}
\hline
\multicolumn{3}{|l|}{\textbf{Level $N = 22$.  Nebentype $\eta = 1$.  Field $\mathbb{F} = GF(16001^{2})$.}} \\
\multicolumn{3}{|l|}{Computed $T_{3}$, $T_{5}$, $T_{7}$.} \\
\hline
1 & $3$ & $\varepsilon^{0} \oplus \varepsilon^{1} \oplus \varepsilon^{2} \sigma_{11,2}$ \\
\hline
1 & $3$ & $\varepsilon^{2} \oplus \varepsilon^{3} \oplus \varepsilon^{0} \sigma_{11,2}$ \\
\hline
1 & $1$ & $\varepsilon^{1} \oplus \varepsilon^{2} \oplus \varepsilon^{0} \sigma_{22,4}$ \\
\hline
\end{tabular}
\end{center}

\begin{center}
\begin{tabular}{|r|r|c|}
\hline
\multicolumn{3}{|l|}{\textbf{Level $N = 22$.  Nebentype $\eta = \chi_{22}^{2}$.  Field $\mathbb{F} = GF(16001^{2})$.}} \\
\multicolumn{3}{|l|}{Computed $T_{3}$, $T_{5}$, $T_{7}$.} \\
\hline
1 & $1$ & $\varepsilon^{0} \oplus \varepsilon^{1} \oplus \varepsilon^{2} \sigma_{22,2}$ \\
\hline
1 & $1$ & $\varepsilon^{2} \oplus \varepsilon^{3} \oplus \varepsilon^{0} \sigma_{22,2}$ \\
\hline
\end{tabular}
\end{center}

\begin{center}
\begin{tabular}{|r|r|c|}
\hline
\multicolumn{3}{|l|}{\textbf{Level $N = 23$.  Nebentype $\eta = 1$.  Field $\mathbb{F} = GF(22067^{60})$.}} \\
\multicolumn{3}{|l|}{Computed $T_{2}$, $T_{3}$, $T_{5}$, $T_{7}$.} \\
\hline
2 & $1$ & $\varepsilon^{0} \oplus \varepsilon^{1} \oplus \varepsilon^{2} \sigma_{23,2a}$ \\
\hline
2 & $1$ & $\varepsilon^{2} \oplus \varepsilon^{3} \oplus \varepsilon^{0} \sigma_{23,2a}$ \\
\hline
1 & $1$ & $\varepsilon^{1} \oplus \varepsilon^{2} \oplus \varepsilon^{0} \sigma_{23,4}$ \\
\hline
\end{tabular}
\end{center}

\begin{center}
\begin{tabular}{|r|r|c|}
\hline
\multicolumn{3}{|l|}{\textbf{Level $N = 23$.  Nebentype $\eta = \chi_{23}^{2}$.  Field $\mathbb{F} = GF(22067^{60})$.}} \\
\multicolumn{3}{|l|}{Computed $T_{2}$, $T_{3}$, $T_{5}$, $T_{7}$.} \\
\hline
1 & $1$ & $\varepsilon^{0} \oplus \varepsilon^{1} \oplus \varepsilon^{2} \sigma_{23,2b}$ \\
\hline
1 & $1$ & $\varepsilon^{2} \oplus \varepsilon^{3} \oplus \varepsilon^{0} \sigma_{23,2b}$ \\
\hline
\end{tabular}
\end{center}

\begin{center}
\begin{tabular}{|r|r|c|}
\hline
\multicolumn{3}{|l|}{\textbf{Level $N = 24$.  Nebentype $\eta = 1$.  Field $\mathbb{F} = GF(12379^{2})$.}} \\
\multicolumn{3}{|l|}{Computed $T_{5}$, $T_{7,1}$, $T_{11,1}$, $T_{13,1}$.} \\
\hline
1 & $1$ & $\varepsilon^{0} \oplus \varepsilon^{1} \oplus \varepsilon^{2} \sigma_{24,2a}$ \\
\hline
1 & $1$ & $\varepsilon^{2} \oplus \varepsilon^{3} \oplus \varepsilon^{0} \sigma_{24,2a}$ \\
\hline
\end{tabular}
\end{center}

\begin{center}
\begin{tabular}{|r|r|c|}
\hline
\multicolumn{3}{|l|}{\textbf{Level $N = 24$.  Nebentype $\eta = \chi_{24,1}$.  Field $\mathbb{F} = GF(12379^{2})$.}} \\
\multicolumn{3}{|l|}{Computed $T_{5}$, $T_{7,1}$, $T_{11,1}$, $T_{13,1}$.} \\
\hline
2 & $1$ & $\varepsilon^{0} \oplus \varepsilon^{1} \oplus \varepsilon^{2} \sigma_{24,2b}$ \\
\hline
2 & $1$ & $\varepsilon^{2} \oplus \varepsilon^{3} \oplus \varepsilon^{0} \sigma_{24,2b}$ \\
\hline
\end{tabular}
\end{center}

\begin{center}
\begin{tabular}{|r|r|c|}
\hline
\multicolumn{3}{|l|}{\textbf{Level $N = 24$.  Nebentype $\eta = \chi_{24,0}\chi_{24,2}$.  Field $\mathbb{F} = GF(12379^{2})$.}} \\
\multicolumn{3}{|l|}{Computed $T_{5}$, $T_{7,1}$, $T_{11,1}$, $T_{13,1}$.} \\
\hline
1 & $3$ & $\varepsilon^{0} \oplus \varepsilon^{1} \oplus \chi_{24,0} \varepsilon^{2} \oplus \chi_{24,2} \varepsilon^{3}$ \\
\hline
1 & $3$ & $\varepsilon^{0} \oplus \varepsilon^{1} \oplus \chi_{24,2} \varepsilon^{2} \oplus \chi_{24,0} \varepsilon^{3}$ \\
\hline
1 & $3$ & $\chi_{24,0} \varepsilon^{0} \oplus \varepsilon^{1} \oplus \varepsilon^{2} \oplus \chi_{24,2} \varepsilon^{3}$ \\
\hline
1 & $3$ & $\chi_{24,2} \varepsilon^{0} \oplus \varepsilon^{1} \oplus \varepsilon^{2} \oplus \chi_{24,0} \varepsilon^{3}$ \\
\hline
1 & $3$ & $\chi_{24,0} \varepsilon^{0} \oplus \chi_{24,2} \varepsilon^{1} \oplus \varepsilon^{2} \oplus \varepsilon^{3}$ \\
\hline
1 & $3$ & $\chi_{24,2} \varepsilon^{0} \oplus \chi_{24,0} \varepsilon^{1} \oplus \varepsilon^{2} \oplus \varepsilon^{3}$ \\
\hline
\end{tabular}
\end{center}

\begin{center}
\begin{tabular}{|r|r|c|}
\hline
\multicolumn{3}{|l|}{\textbf{Level $N = 24$.  Nebentype $\eta = \chi_{24,0}\chi_{24,1}\chi_{24,2}$.  Field $\mathbb{F} = GF(12379^{2})$.}} \\
\multicolumn{3}{|l|}{Computed $T_{5}$, $T_{7,1}$, $T_{11,1}$, $T_{13,1}$, $T_{17,1}$.} \\
\hline
1 & $1$ & $\varepsilon^{0} \oplus \varepsilon^{1} \oplus \chi_{24,0}\chi_{24,1} \varepsilon^{2} \oplus \chi_{24,2} \varepsilon^{3}$ \\
\hline
1 & $1$ & $\varepsilon^{0} \oplus \varepsilon^{1} \oplus \chi_{24,2} \varepsilon^{2} \oplus \chi_{24,0}\chi_{24,1} \varepsilon^{3}$ \\
\hline
1 & $1$ & $\chi_{24,0}\chi_{24,1} \varepsilon^{0} \oplus \varepsilon^{1} \oplus \varepsilon^{2} \oplus \chi_{24,2} \varepsilon^{3}$ \\
\hline
1 & $1$ & $\chi_{24,2} \varepsilon^{0} \oplus \varepsilon^{1} \oplus \varepsilon^{2} \oplus \chi_{24,0}\chi_{24,1} \varepsilon^{3}$ \\
\hline
1 & $1$ & $\chi_{24,0}\chi_{24,1} \varepsilon^{0} \oplus \chi_{24,2} \varepsilon^{1} \oplus \varepsilon^{2} \oplus \varepsilon^{3}$ \\
\hline
1 & $1$ & $\chi_{24,2} \varepsilon^{0} \oplus \chi_{24,0}\chi_{24,1} \varepsilon^{1} \oplus \varepsilon^{2} \oplus \varepsilon^{3}$ \\
\hline
2 & $1$ & $\varepsilon^{0} \oplus \varepsilon^{1} \oplus \varepsilon^{2} \sigma_{24,2c}$ \\
\hline
2 & $1$ & $\varepsilon^{2} \oplus \varepsilon^{3} \oplus \varepsilon^{0} \sigma_{24,2c}$ \\
\hline
1 & $1$ & $\varepsilon^{0} \oplus \chi_{24,2} \varepsilon^{2} \oplus \varepsilon^{1} \sigma_{8,3}$ \\
\hline
1 & $1$ & $\varepsilon^{1} \oplus \chi_{24,2} \varepsilon^{3} \oplus \varepsilon^{0} \sigma_{8,3}$ \\
\hline
1 & $1$ & $\chi_{24,2} \varepsilon^{0} \oplus \varepsilon^{2} \oplus \varepsilon^{1} \sigma_{8,3}$ \\
\hline
1 & $1$ & $\chi_{24,2} \varepsilon^{1} \oplus \varepsilon^{3} \oplus \varepsilon^{0} \sigma_{8,3}$ \\
\hline
\end{tabular}
\end{center}

\begin{center}
\begin{tabular}{|r|r|c|}
\hline
\multicolumn{3}{|l|}{\textbf{Level $N = 25$.  Nebentype $\eta = 1$.  Field $\mathbb{F} = GF(16001^{60})$.}} \\
\multicolumn{3}{|l|}{Computed $T_{2}$, $T_{3}$.} \\
\hline
1 & $1$ & $\varepsilon^{0} \oplus \varepsilon^{1} \oplus \chi_{25}^{15} \varepsilon^{2} \oplus \chi_{25}^{5} \varepsilon^{3}$ \\
\hline
1 & $1$ & $\varepsilon^{0} \oplus \varepsilon^{1} \oplus \chi_{25}^{5} \varepsilon^{2} \oplus \chi_{25}^{15} \varepsilon^{3}$ \\
\hline
1 & $1$ & $\chi_{25}^{15} \varepsilon^{0} \oplus \varepsilon^{1} \oplus \varepsilon^{2} \oplus \chi_{25}^{5} \varepsilon^{3}$ \\
\hline
1 & $1$ & $\chi_{25}^{5} \varepsilon^{0} \oplus \varepsilon^{1} \oplus \varepsilon^{2} \oplus \chi_{25}^{15} \varepsilon^{3}$ \\
\hline
1 & $1$ & $\chi_{25}^{15} \varepsilon^{0} \oplus \chi_{25}^{5} \varepsilon^{1} \oplus \varepsilon^{2} \oplus \varepsilon^{3}$ \\
\hline
1 & $1$ & $\chi_{25}^{5} \varepsilon^{0} \oplus \chi_{25}^{15} \varepsilon^{1} \oplus \varepsilon^{2} \oplus \varepsilon^{3}$ \\
\hline
1 & $1$ & $\varepsilon^{1} \oplus \varepsilon^{2} \oplus \varepsilon^{0} \sigma_{25,4}$ \\
\hline
\end{tabular}
\end{center}

\begin{center}
\begin{tabular}{|r|r|c|}
\hline
\multicolumn{3}{|l|}{\textbf{Level $N = 25$.  Nebentype $\eta = \chi_{25}^{2}$.  Field $\mathbb{F} = GF(16001^{60})$.}} \\
\multicolumn{3}{|l|}{Computed $T_{2}$, $T_{3}$.} \\
\hline
2 & $1$ & $\varepsilon^{0} \oplus \varepsilon^{1} \oplus \varepsilon^{2} \sigma_{25,2a}$ \\
\hline
2 & $1$ & $\varepsilon^{2} \oplus \varepsilon^{3} \oplus \varepsilon^{0} \sigma_{25,2a}$ \\
\hline
\end{tabular}
\end{center}

\begin{center}
\begin{tabular}{|r|r|c|}
\hline
\multicolumn{3}{|l|}{\textbf{Level $N = 25$.  Nebentype $\eta = \chi_{25}^{4}$.  Field $\mathbb{F} = GF(16001^{60})$.}} \\
\multicolumn{3}{|l|}{Computed $T_{2}$, $T_{3}$.} \\
\hline
1 & $1$ & $\varepsilon^{0} \oplus \varepsilon^{1} \oplus \varepsilon^{2} \sigma_{25,2b}$ \\
\hline
1 & $1$ & $\varepsilon^{2} \oplus \varepsilon^{3} \oplus \varepsilon^{0} \sigma_{25,2b}$ \\
\hline
\end{tabular}
\end{center}

\begin{center}
\begin{tabular}{|r|r|c|}
\hline
\multicolumn{3}{|l|}{\textbf{Level $N = 25$.  Nebentype $\eta = \chi_{25}^{10}$.  Field $\mathbb{F} = GF(16001^{60})$.}} \\
\multicolumn{3}{|l|}{Computed $T_{2}$, $T_{3}$.} \\
\hline
1 & $1$ & $\varepsilon^{0} \oplus \varepsilon^{1} \oplus \chi_{25}^{15} \varepsilon^{2} \oplus \chi_{25}^{15} \varepsilon^{3}$ \\
\hline
1 & $1$ & $\varepsilon^{0} \oplus \varepsilon^{1} \oplus \chi_{25}^{5} \varepsilon^{2} \oplus \chi_{25}^{5} \varepsilon^{3}$ \\
\hline
1 & $1$ & $\chi_{25}^{15} \varepsilon^{0} \oplus \varepsilon^{1} \oplus \varepsilon^{2} \oplus \chi_{25}^{15} \varepsilon^{3}$ \\
\hline
1 & $1$ & $\chi_{25}^{5} \varepsilon^{0} \oplus \varepsilon^{1} \oplus \varepsilon^{2} \oplus \chi_{25}^{5} \varepsilon^{3}$ \\
\hline
1 & $1$ & $\chi_{25}^{15} \varepsilon^{0} \oplus \chi_{25}^{15} \varepsilon^{1} \oplus \varepsilon^{2} \oplus \varepsilon^{3}$ \\
\hline
1 & $1$ & $\chi_{25}^{5} \varepsilon^{0} \oplus \chi_{25}^{5} \varepsilon^{1} \oplus \varepsilon^{2} \oplus \varepsilon^{3}$ \\
\hline
\end{tabular}
\end{center}

\begin{center}
\begin{tabular}{|r|r|c|}
\hline
\multicolumn{3}{|l|}{\textbf{Level $N = 26$.  Nebentype $\eta = 1$.  Field $\mathbb{F} = GF(12037^{2})$.}} \\
\multicolumn{3}{|l|}{Computed $T_{3}$, $T_{5}$.} \\
\hline
1 & $1$ & $\varepsilon^{0} \oplus \varepsilon^{1} \oplus \varepsilon^{2} \sigma_{26,2a}$ \\
\hline
1 & $1$ & $\varepsilon^{2} \oplus \varepsilon^{3} \oplus \varepsilon^{0} \sigma_{26,2a}$ \\
\hline
1 & $1$ & $\varepsilon^{0} \oplus \varepsilon^{1} \oplus \varepsilon^{2} \sigma_{26,2b}$ \\
\hline
1 & $1$ & $\varepsilon^{2} \oplus \varepsilon^{3} \oplus \varepsilon^{0} \sigma_{26,2b}$ \\
\hline
1 & $3$ & $\varepsilon^{1} \oplus \varepsilon^{2} \oplus \varepsilon^{0} \sigma_{13,4}$ \\
\hline
\end{tabular}
\end{center}

\begin{center}
\begin{tabular}{|r|r|c|}
\hline
\multicolumn{3}{|l|}{\textbf{Level $N = 26$.  Nebentype $\eta = \chi_{26}^{2}$.  Field $\mathbb{F} = GF(12037^{2})$.}} \\
\multicolumn{3}{|l|}{Computed $T_{3}$, $T_{5}$.} \\
\hline
1 & $3$ & $\varepsilon^{0} \oplus \varepsilon^{1} \oplus \varepsilon^{2} \sigma_{13,2}$ \\
\hline
1 & $3$ & $\varepsilon^{2} \oplus \varepsilon^{3} \oplus \varepsilon^{0} \sigma_{13,2}$ \\
\hline
\end{tabular}
\end{center}

\begin{center}
\begin{tabular}{|r|r|c|}
\hline
\multicolumn{3}{|l|}{\textbf{Level $N = 26$.  Nebentype $\eta = \chi_{26}^{4}$.  Field $\mathbb{F} = GF(12037^{2})$.}} \\
\multicolumn{3}{|l|}{Computed $T_{3}$, $T_{5}$.} \\
\hline
1 & $1$ & $\varepsilon^{0} \oplus \varepsilon^{1} \oplus \varepsilon^{2} \sigma_{26,2c}$ \\
\hline
1 & $1$ & $\varepsilon^{2} \oplus \varepsilon^{3} \oplus \varepsilon^{0} \sigma_{26,2c}$ \\
\hline
\end{tabular}
\end{center}

\begin{center}
\begin{tabular}{|r|r|c|}
\hline
\multicolumn{3}{|l|}{\textbf{Level $N = 26$.  Nebentype $\eta = \chi_{26}^{6}$.  Field $\mathbb{F} = GF(12037^{2})$.}} \\
\multicolumn{3}{|l|}{Computed $T_{3}$, $T_{5}$.} \\
\hline
2 & $1$ & $\varepsilon^{0} \oplus \varepsilon^{1} \oplus \varepsilon^{2} \sigma_{26,2d}$ \\
\hline
2 & $1$ & $\varepsilon^{2} \oplus \varepsilon^{3} \oplus \varepsilon^{0} \sigma_{26,2d}$ \\
\hline
\end{tabular}
\end{center}

\begin{center}
\begin{tabular}{|r|r|c|}
\hline
\multicolumn{3}{|l|}{\textbf{Level $N = 27$.  Nebentype $\eta = 1$.  Field $\mathbb{F} = GF(11863^{6})$.}} \\
\multicolumn{3}{|l|}{Computed $T_{2}$, $T_{5}$, $T_{7,1}$.} \\
\hline
1 & $3$ & $\varepsilon^{0} \oplus \varepsilon^{1} \oplus \chi_{27}^{9} \varepsilon^{2} \oplus \chi_{27}^{9} \varepsilon^{3}$ \\
\hline
1 & $3$ & $\chi_{27}^{9} \varepsilon^{0} \oplus \varepsilon^{1} \oplus \varepsilon^{2} \oplus \chi_{27}^{9} \varepsilon^{3}$ \\
\hline
1 & $3$ & $\chi_{27}^{9} \varepsilon^{0} \oplus \chi_{27}^{9} \varepsilon^{1} \oplus \varepsilon^{2} \oplus \varepsilon^{3}$ \\
\hline
1 & $1$ & $\varepsilon^{0} \oplus \varepsilon^{1} \oplus \varepsilon^{2} \sigma_{27,2a}$ \\
\hline
1 & $1$ & $\varepsilon^{2} \oplus \varepsilon^{3} \oplus \varepsilon^{0} \sigma_{27,2a}$ \\
\hline
1 & $1$ & $\varepsilon^{1} \oplus \varepsilon^{2} \oplus \varepsilon^{0} \sigma_{27,4}$ \\
\hline
\end{tabular}
\end{center}

\begin{center}
\begin{tabular}{|r|r|c|}
\hline
\multicolumn{3}{|l|}{\textbf{Level $N = 27$.  Nebentype $\eta = \chi_{27}^{2}$.  Field $\mathbb{F} = GF(11863^{6})$.}} \\
\multicolumn{3}{|l|}{Computed $T_{2}$, $T_{5}$.} \\
\hline
2 & $1$ & $\varepsilon^{0} \oplus \varepsilon^{1} \oplus \varepsilon^{2} \sigma_{27,2b}$ \\
\hline
2 & $1$ & $\varepsilon^{2} \oplus \varepsilon^{3} \oplus \varepsilon^{0} \sigma_{27,2b}$ \\
\hline
\end{tabular}
\end{center}

\begin{center}
\begin{tabular}{|r|r|c|}
\hline
\multicolumn{3}{|l|}{\textbf{Level $N = 27$.  Nebentype $\eta = \chi_{27}^{6}$.  Field $\mathbb{F} = GF(11863^{6})$.}} \\
\multicolumn{3}{|l|}{Computed $T_{2}$, $T_{5}$, $T_{7,1}$.} \\
\hline
1 & $1$ & $\varepsilon^{0} \oplus \varepsilon^{1} \oplus \chi_{27}^{15} \varepsilon^{2} \oplus \chi_{27}^{9} \varepsilon^{3}$ \\
\hline
1 & $1$ & $\varepsilon^{0} \oplus \varepsilon^{1} \oplus \chi_{27}^{9} \varepsilon^{2} \oplus \chi_{27}^{15} \varepsilon^{3}$ \\
\hline
1 & $1$ & $\chi_{27}^{15} \varepsilon^{0} \oplus \varepsilon^{1} \oplus \varepsilon^{2} \oplus \chi_{27}^{9} \varepsilon^{3}$ \\
\hline
1 & $1$ & $\chi_{27}^{9} \varepsilon^{0} \oplus \varepsilon^{1} \oplus \varepsilon^{2} \oplus \chi_{27}^{15} \varepsilon^{3}$ \\
\hline
1 & $1$ & $\chi_{27}^{15} \varepsilon^{0} \oplus \chi_{27}^{9} \varepsilon^{1} \oplus \varepsilon^{2} \oplus \varepsilon^{3}$ \\
\hline
1 & $1$ & $\chi_{27}^{9} \varepsilon^{0} \oplus \chi_{27}^{15} \varepsilon^{1} \oplus \varepsilon^{2} \oplus \varepsilon^{3}$ \\
\hline
1 & $1$ & $\varepsilon^{0} \oplus \chi_{27}^{9} \varepsilon^{2} \oplus \varepsilon^{1} \sigma_{9,3}$ \\
\hline
1 & $1$ & $\varepsilon^{1} \oplus \chi_{27}^{9} \varepsilon^{3} \oplus \varepsilon^{0} \sigma_{9,3}$ \\
\hline
1 & $1$ & $\chi_{27}^{9} \varepsilon^{0} \oplus \varepsilon^{2} \oplus \varepsilon^{1} \sigma_{9,3}$ \\
\hline
1 & $1$ & $\chi_{27}^{9} \varepsilon^{1} \oplus \varepsilon^{3} \oplus \varepsilon^{0} \sigma_{9,3}$ \\
\hline
\end{tabular}
\end{center}

\begin{center}
\begin{tabular}{|r|r|c|}
\hline
\multicolumn{3}{|l|}{\textbf{Level $N = 28$.  Nebentype $\eta = 1$.  Field $\mathbb{F} = GF(12379^{12})$.}} \\
\multicolumn{3}{|l|}{Computed $T_{3}$, $T_{5}$, $T_{11,1}$, $T_{13,1}$.} \\
\hline
1 & $3$ & $\varepsilon^{0} \oplus \varepsilon^{1} \oplus \varepsilon^{2} \sigma_{14,2}$ \\
\hline
1 & $3$ & $\varepsilon^{2} \oplus \varepsilon^{3} \oplus \varepsilon^{0} \sigma_{14,2}$ \\
\hline
1 & $1$ & $\varepsilon^{1} \oplus \varepsilon^{2} \oplus \varepsilon^{0} \sigma_{28,4}$ \\
\hline
\end{tabular}
\end{center}

\begin{center}
\begin{tabular}{|r|r|c|}
\hline
\multicolumn{3}{|l|}{\textbf{Level $N = 28$.  Nebentype $\eta = \chi_{28,1}^{2}$.  Field $\mathbb{F} = GF(12379^{12})$.}} \\
\multicolumn{3}{|l|}{Computed $T_{3}$, $T_{5}$.} \\
\hline
1 & $1$ & $\varepsilon^{0} \oplus \varepsilon^{1} \oplus \varepsilon^{2} \sigma_{28,2a}$ \\
\hline
1 & $1$ & $\varepsilon^{2} \oplus \varepsilon^{3} \oplus \varepsilon^{0} \sigma_{28,2a}$ \\
\hline
\end{tabular}
\end{center}

\begin{center}
\begin{tabular}{|r|r|c|}
\hline
\multicolumn{3}{|l|}{\textbf{Level $N = 28$.  Nebentype $\eta = \chi_{28,0}\chi_{28,1}$.  Field $\mathbb{F} = GF(12379^{12})$.}} \\
\multicolumn{3}{|l|}{Computed $T_{3}$, $T_{5}$.} \\
\hline
1 & $1$ & $\varepsilon^{0} \oplus \varepsilon^{1} \oplus \chi_{28,0} \varepsilon^{2} \oplus \chi_{28,1} \varepsilon^{3}$ \\
\hline
1 & $1$ & $\varepsilon^{0} \oplus \varepsilon^{1} \oplus \chi_{28,1} \varepsilon^{2} \oplus \chi_{28,0} \varepsilon^{3}$ \\
\hline
1 & $1$ & $\chi_{28,0} \varepsilon^{0} \oplus \varepsilon^{1} \oplus \varepsilon^{2} \oplus \chi_{28,1} \varepsilon^{3}$ \\
\hline
1 & $1$ & $\chi_{28,1} \varepsilon^{0} \oplus \varepsilon^{1} \oplus \varepsilon^{2} \oplus \chi_{28,0} \varepsilon^{3}$ \\
\hline
1 & $1$ & $\chi_{28,0} \varepsilon^{0} \oplus \chi_{28,1} \varepsilon^{1} \oplus \varepsilon^{2} \oplus \varepsilon^{3}$ \\
\hline
1 & $1$ & $\chi_{28,1} \varepsilon^{0} \oplus \chi_{28,0} \varepsilon^{1} \oplus \varepsilon^{2} \oplus \varepsilon^{3}$ \\
\hline
2 & $1$ & $\varepsilon^{0} \oplus \varepsilon^{1} \oplus \varepsilon^{2} \sigma_{28,2b}$ \\
\hline
2 & $1$ & $\varepsilon^{2} \oplus \varepsilon^{3} \oplus \varepsilon^{0} \sigma_{28,2b}$ \\
\hline
\end{tabular}
\end{center}

\begin{center}
\begin{tabular}{|r|r|c|}
\hline
\multicolumn{3}{|l|}{\textbf{Level $N = 28$.  Nebentype $\eta = \chi_{28,0}\chi_{28,1}^{3}$.  Field $\mathbb{F} = GF(12379^{12})$.}} \\
\multicolumn{3}{|l|}{Computed $T_{3}$, $T_{5}$, $T_{11,1}$, $T_{13,1}$.} \\
\hline
1 & $1$ & $\varepsilon^{0} \oplus \varepsilon^{1} \oplus \chi_{28,0} \varepsilon^{2} \oplus \chi_{28,1}^{3} \varepsilon^{3}$ \\
\hline
1 & $1$ & $\varepsilon^{0} \oplus \varepsilon^{1} \oplus \chi_{28,1}^{3} \varepsilon^{2} \oplus \chi_{28,0} \varepsilon^{3}$ \\
\hline
1 & $1$ & $\chi_{28,0} \varepsilon^{0} \oplus \varepsilon^{1} \oplus \varepsilon^{2} \oplus \chi_{28,1}^{3} \varepsilon^{3}$ \\
\hline
1 & $1$ & $\chi_{28,1}^{3} \varepsilon^{0} \oplus \varepsilon^{1} \oplus \varepsilon^{2} \oplus \chi_{28,0} \varepsilon^{3}$ \\
\hline
1 & $1$ & $\chi_{28,0} \varepsilon^{0} \oplus \chi_{28,1}^{3} \varepsilon^{1} \oplus \varepsilon^{2} \oplus \varepsilon^{3}$ \\
\hline
1 & $1$ & $\chi_{28,1}^{3} \varepsilon^{0} \oplus \chi_{28,0} \varepsilon^{1} \oplus \varepsilon^{2} \oplus \varepsilon^{3}$ \\
\hline
2 & $1$ & $\varepsilon^{0} \oplus \varepsilon^{1} \oplus \varepsilon^{2} \sigma_{28,2c}$ \\
\hline
2 & $1$ & $\varepsilon^{2} \oplus \varepsilon^{3} \oplus \varepsilon^{0} \sigma_{28,2c}$ \\
\hline
1 & $1$ & $\varepsilon^{0} \oplus \chi_{28,0} \varepsilon^{2} \oplus \varepsilon^{1} \sigma_{7,3}$ \\
\hline
1 & $1$ & $\varepsilon^{1} \oplus \chi_{28,0} \varepsilon^{3} \oplus \varepsilon^{0} \sigma_{7,3}$ \\
\hline
1 & $1$ & $\chi_{28,0} \varepsilon^{0} \oplus \varepsilon^{2} \oplus \varepsilon^{1} \sigma_{7,3}$ \\
\hline
1 & $1$ & $\chi_{28,0} \varepsilon^{1} \oplus \varepsilon^{3} \oplus \varepsilon^{0} \sigma_{7,3}$ \\
\hline
\end{tabular}
\end{center}

\begin{center}
\begin{tabular}{|r|r|c|}
\hline
\multicolumn{3}{|l|}{\textbf{Level $N = 29$.  Nebentype $\eta = 1$.  Field $\mathbb{F} = GF(2297^{6})$.}} \\
\multicolumn{3}{|l|}{Computed $T_{2}$, $T_{3}$, $T_{5}$.} \\
\hline
2 & $1$ & $\varepsilon^{0} \oplus \varepsilon^{1} \oplus \varepsilon^{2} \sigma_{29,2a}$ \\
\hline
2 & $1$ & $\varepsilon^{2} \oplus \varepsilon^{3} \oplus \varepsilon^{0} \sigma_{29,2a}$ \\
\hline
2 & $1$ & $\varepsilon^{1} \oplus \varepsilon^{2} \oplus \varepsilon^{0} \sigma_{29,4}$ \\
\hline
\end{tabular}
\end{center}

\begin{center}
\begin{tabular}{|r|r|c|}
\hline
\multicolumn{3}{|l|}{\textbf{Level $N = 29$.  Nebentype $\eta = \chi_{29}^{2}$.  Field $\mathbb{F} = GF(2297^{6})$.}} \\
\multicolumn{3}{|l|}{Computed $T_{2}$, $T_{3}$, $T_{5}$.} \\
\hline
2 & $1$ & $\varepsilon^{0} \oplus \varepsilon^{1} \oplus \varepsilon^{2} \sigma_{29,2b}$ \\
\hline
2 & $1$ & $\varepsilon^{2} \oplus \varepsilon^{3} \oplus \varepsilon^{0} \sigma_{29,2b}$ \\
\hline
\end{tabular}
\end{center}

\begin{center}
\begin{tabular}{|r|r|c|}
\hline
\multicolumn{3}{|l|}{\textbf{Level $N = 29$.  Nebentype $\eta = \chi_{29}^{4}$.  Field $\mathbb{F} = GF(2297^{6})$.}} \\
\multicolumn{3}{|l|}{Computed $T_{2}$, $T_{3}$, $T_{5}$.} \\
\hline
1 & $1$ & $\varepsilon^{0} \oplus \varepsilon^{1} \oplus \varepsilon^{2} \sigma_{29,2c}$ \\
\hline
1 & $1$ & $\varepsilon^{2} \oplus \varepsilon^{3} \oplus \varepsilon^{0} \sigma_{29,2c}$ \\
\hline
\end{tabular}
\end{center}

\begin{center}
\begin{tabular}{|r|r|c|}
\hline
\multicolumn{3}{|l|}{\textbf{Level $N = 29$.  Nebentype $\eta = \chi_{29}^{14}$.  Field $\mathbb{F} = GF(2297^{6})$.}} \\
\multicolumn{3}{|l|}{Computed $T_{2}$, $T_{3}$, $T_{5}$.} \\
\hline
2 & $1$ & $\varepsilon^{0} \oplus \varepsilon^{1} \oplus \varepsilon^{2} \sigma_{29,2d}$ \\
\hline
2 & $1$ & $\varepsilon^{2} \oplus \varepsilon^{3} \oplus \varepsilon^{0} \sigma_{29,2d}$ \\
\hline
2 & $1$ & $\varepsilon^{0} \oplus \varepsilon^{1} \mathrm{Sym}^2(\sigma_{29,2d})$ \\
\hline
\end{tabular}
\end{center}

\begin{center}
\begin{tabular}{|r|r|c|}
\hline
\multicolumn{3}{|l|}{\textbf{Level $N = 31$.  Nebentype $\eta = 1$.  Field $\mathbb{F} = GF(4201^{60})$.}} \\
\multicolumn{3}{|l|}{Computed $T_{2}$, $T_{3}$, $T_{5}$.} \\
\hline
2 & $1$ & $\varepsilon^{0} \oplus \varepsilon^{1} \oplus \varepsilon^{2} \sigma_{31,2a}$ \\
\hline
2 & $1$ & $\varepsilon^{2} \oplus \varepsilon^{3} \oplus \varepsilon^{0} \sigma_{31,2a}$ \\
\hline
2 & $1$ & $\varepsilon^{1} \oplus \varepsilon^{2} \oplus \varepsilon^{0} \sigma_{31,4}$ \\
\hline
\end{tabular}
\end{center}

\begin{center}
\begin{tabular}{|r|r|c|}
\hline
\multicolumn{3}{|l|}{\textbf{Level $N = 31$.  Nebentype $\eta = \chi_{31}^{2}$.  Field $\mathbb{F} = GF(4201^{60})$.}} \\
\multicolumn{3}{|l|}{Computed $T_{2}$, $T_{3}$, $T_{5}$.} \\
\hline
2 & $1$ & $\varepsilon^{0} \oplus \varepsilon^{1} \oplus \varepsilon^{2} \sigma_{31,2b}$ \\
\hline
2 & $1$ & $\varepsilon^{2} \oplus \varepsilon^{3} \oplus \varepsilon^{0} \sigma_{31,2b}$ \\
\hline
\end{tabular}
\end{center}

\begin{center}
\begin{tabular}{|r|r|c|}
\hline
\multicolumn{3}{|l|}{\textbf{Level $N = 31$.  Nebentype $\eta = \chi_{31}^{6}$.  Field $\mathbb{F} = GF(4201^{60})$.}} \\
\multicolumn{3}{|l|}{Computed $T_{2}$, $T_{3}$, $T_{5}$.} \\
\hline
1 & $1$ & $\varepsilon^{0} \oplus \varepsilon^{1} \oplus \varepsilon^{2} \sigma_{31,2c}$ \\
\hline
1 & $1$ & $\varepsilon^{2} \oplus \varepsilon^{3} \oplus \varepsilon^{0} \sigma_{31,2c}$ \\
\hline
\end{tabular}
\end{center}

\begin{center}
\begin{tabular}{|r|r|c|}
\hline
\multicolumn{3}{|l|}{\textbf{Level $N = 31$.  Nebentype $\eta = \chi_{31}^{10}$.  Field $\mathbb{F} = GF(4201^{60})$.}} \\
\multicolumn{3}{|l|}{Computed $T_{2}$, $T_{3}$, $T_{5}$.} \\
\hline
2 & $1$ & $\varepsilon^{0} \oplus \varepsilon^{1} \oplus \varepsilon^{2} \sigma_{31,2d}$ \\
\hline
2 & $1$ & $\varepsilon^{2} \oplus \varepsilon^{3} \oplus \varepsilon^{0} \sigma_{31,2d}$ \\
\hline
\end{tabular}
\end{center}

\begin{center}
\begin{tabular}{|r|r|c|}
\hline
\multicolumn{3}{|l|}{\textbf{Level $N = 37$.  Nebentype $\eta = 1$.  Field $\mathbb{F} = GF(3889^{24})$.}} \\
\multicolumn{3}{|l|}{Computed $T_{2}$, $T_{3}$, $T_{5}$, $T_{7,1}$, $T_{13,1}$.} \\
\hline
1 & $1$ & $\varepsilon^{0} \oplus \varepsilon^{1} \oplus \varepsilon^{2} \sigma_{37,2a}$ \\
\hline
1 & $1$ & $\varepsilon^{2} \oplus \varepsilon^{3} \oplus \varepsilon^{0} \sigma_{37,2a}$ \\
\hline
1 & $1$ & $\varepsilon^{0} \oplus \varepsilon^{1} \oplus \varepsilon^{2} \sigma_{37,2b}$ \\
\hline
1 & $1$ & $\varepsilon^{2} \oplus \varepsilon^{3} \oplus \varepsilon^{0} \sigma_{37,2b}$ \\
\hline
4 & $1$ & $\varepsilon^{1} \oplus \varepsilon^{2} \oplus \varepsilon^{0} \sigma_{37,4}$ \\
\hline
\end{tabular}
\end{center}

\begin{center}
\begin{tabular}{|r|r|c|}
\hline
\multicolumn{3}{|l|}{\textbf{Level $N = 37$.  Nebentype $\eta = \chi_{37}^{2}$.  Field $\mathbb{F} = GF(3889^{24})$.}} \\
\multicolumn{3}{|l|}{Computed $T_{2}$, $T_{3}$, $T_{5}$.} \\
\hline
3 & $1$ & $\varepsilon^{0} \oplus \varepsilon^{1} \oplus \varepsilon^{2} \sigma_{37,2c}$ \\
\hline
3 & $1$ & $\varepsilon^{2} \oplus \varepsilon^{3} \oplus \varepsilon^{0} \sigma_{37,2c}$ \\
\hline
\end{tabular}
\end{center}

\begin{center}
\begin{tabular}{|r|r|c|}
\hline
\multicolumn{3}{|l|}{\textbf{Level $N = 37$.  Nebentype $\eta = \chi_{37}^{4}$.  Field $\mathbb{F} = GF(3889^{24})$.}} \\
\multicolumn{3}{|l|}{Computed $T_{2}$, $T_{3}$, $T_{5}$.} \\
\hline
1 & $1$ & $\varepsilon^{0} \oplus \varepsilon^{1} \oplus \varepsilon^{2} \sigma_{37,2d}$ \\
\hline
1 & $1$ & $\varepsilon^{2} \oplus \varepsilon^{3} \oplus \varepsilon^{0} \sigma_{37,2d}$ \\
\hline
1 & $1$ & $\varepsilon^{0} \oplus \varepsilon^{1} \oplus \varepsilon^{2} \sigma_{37,2e}$ \\
\hline
1 & $1$ & $\varepsilon^{2} \oplus \varepsilon^{3} \oplus \varepsilon^{0} \sigma_{37,2e}$ \\
\hline
\end{tabular}
\end{center}

\begin{center}
\begin{tabular}{|r|r|c|}
\hline
\multicolumn{3}{|l|}{\textbf{Level $N = 37$.  Nebentype $\eta = \chi_{37}^{6}$.  Field $\mathbb{F} = GF(3889^{24})$.}} \\
\multicolumn{3}{|l|}{Computed $T_{2}$, $T_{3}$, $T_{5}$.} \\
\hline
2 & $1$ & $\varepsilon^{0} \oplus \varepsilon^{1} \oplus \varepsilon^{2} \sigma_{37,2{f}}$ \\
\hline
2 & $1$ & $\varepsilon^{2} \oplus \varepsilon^{3} \oplus \varepsilon^{0} \sigma_{37,2{f}}$ \\
\hline
\end{tabular}
\end{center}

\begin{center}
\begin{tabular}{|r|r|c|}
\hline
\multicolumn{3}{|l|}{\textbf{Level $N = 37$.  Nebentype $\eta = \chi_{37}^{12}$.  Field $\mathbb{F} = GF(3889^{24})$.}} \\
\multicolumn{3}{|l|}{Computed $T_{2}$, $T_{3}$, $T_{5}$.} \\
\hline
1 & $1$ & $\varepsilon^{0} \oplus \varepsilon^{1} \oplus \varepsilon^{2} \sigma_{37,2g}$ \\
\hline
1 & $1$ & $\varepsilon^{2} \oplus \varepsilon^{3} \oplus \varepsilon^{0} \sigma_{37,2g}$ \\
\hline
\end{tabular}
\end{center}

\begin{center}
\begin{tabular}{|r|r|c|}
\hline
\multicolumn{3}{|l|}{\textbf{Level $N = 37$.  Nebentype $\eta = \chi_{37}^{18,2}$.  Field $\mathbb{F} = GF(3889^{24})$.}} \\
\multicolumn{3}{|l|}{Computed $T_{2}$, $T_{3}$, $T_{5}$.} \\
\hline
2 & $1$ & $\varepsilon^{0} \oplus \varepsilon^{1} \oplus \varepsilon^{2} \sigma_{37,2h}$ \\
\hline
2 & $1$ & $\varepsilon^{2} \oplus \varepsilon^{3} \oplus \varepsilon^{0} \sigma_{37,2h}$ \\
\hline
2 & $1$ & $\varepsilon^{0} \oplus \varepsilon^{1} \mathrm{Sym}^2(\sigma_{37,2h})$ \\
\hline
\end{tabular}
\end{center}

\begin{center}
\begin{tabular}{|r|r|c|}
\hline
\multicolumn{3}{|l|}{\textbf{Level $N = 41$.  Nebentype $\eta = 1$.  Field $\mathbb{F} = GF(21881^{60})$.}} \\
\multicolumn{3}{|l|}{Computed $T_{2}$, $T_{3}$, $T_{5}$.} \\
\hline
3 & $1$ & $\varepsilon^{0} \oplus \varepsilon^{1} \oplus \varepsilon^{2} \sigma_{41,2a}$ \\
\hline
3 & $1$ & $\varepsilon^{2} \oplus \varepsilon^{3} \oplus \varepsilon^{0} \sigma_{41,2a}$ \\
\hline
3 & $1$ & $\varepsilon^{1} \oplus \varepsilon^{2} \oplus \varepsilon^{0} \sigma_{41,4}$ \\
\hline
\end{tabular}
\end{center}

\begin{center}
\begin{tabular}{|r|r|c|}
\hline
\multicolumn{3}{|l|}{\textbf{Level $N = 41$.  Nebentype $\eta = \chi_{41}^{2}$.  Field $\mathbb{F} = GF(21881^{60})$.}} \\
\multicolumn{3}{|l|}{Computed $T_{2}$, $T_{3}$, $T_{5}$.} \\
\hline
3 & $1$ & $\varepsilon^{0} \oplus \varepsilon^{1} \oplus \varepsilon^{2} \sigma_{41,2b}$ \\
\hline
3 & $1$ & $\varepsilon^{2} \oplus \varepsilon^{3} \oplus \varepsilon^{0} \sigma_{41,2b}$ \\
\hline
\end{tabular}
\end{center}

\begin{center}
\begin{tabular}{|r|r|c|}
\hline
\multicolumn{3}{|l|}{\textbf{Level $N = 41$.  Nebentype $\eta = \chi_{41}^{4}$.  Field $\mathbb{F} = GF(21881^{60})$.}} \\
\multicolumn{3}{|l|}{Computed $T_{2}$, $T_{3}$, $T_{5}$.} \\
\hline
2 & $1$ & $\varepsilon^{0} \oplus \varepsilon^{1} \oplus \varepsilon^{2} \sigma_{41,2c}$ \\
\hline
2 & $1$ & $\varepsilon^{2} \oplus \varepsilon^{3} \oplus \varepsilon^{0} \sigma_{41,2c}$ \\
\hline
\end{tabular}
\end{center}

\begin{center}
\begin{tabular}{|r|r|c|}
\hline
\multicolumn{3}{|l|}{\textbf{Level $N = 41$.  Nebentype $\eta = \chi_{41}^{8}$.  Field $\mathbb{F} = GF(21881^{60})$.}} \\
\multicolumn{3}{|l|}{Computed $T_{2}$, $T_{3}$, $T_{5}$.} \\
\hline
2 & $1$ & $\varepsilon^{0} \oplus \varepsilon^{1} \oplus \varepsilon^{2} \sigma_{41,2d}$ \\
\hline
2 & $1$ & $\varepsilon^{2} \oplus \varepsilon^{3} \oplus \varepsilon^{0} \sigma_{41,2d}$ \\
\hline
\end{tabular}
\end{center}

\begin{center}
\begin{tabular}{|r|r|c|}
\hline
\multicolumn{3}{|l|}{\textbf{Level $N = 41$.  Nebentype $\eta = \chi_{41}^{10}$.  Field $\mathbb{F} = GF(21881^{60})$.}} \\
\multicolumn{3}{|l|}{Computed $T_{2}$, $T_{3}$, $T_{5}$.} \\
\hline
3 & $1$ & $\varepsilon^{0} \oplus \varepsilon^{1} \oplus \varepsilon^{2} \sigma_{41,2e}$ \\
\hline
3 & $1$ & $\varepsilon^{2} \oplus \varepsilon^{3} \oplus \varepsilon^{0} \sigma_{41,2e}$ \\
\hline
1 & $1$ & $\varepsilon^{0} \oplus \varepsilon^{1} \delta$ \\
\hline
1 & $1$ & $\varepsilon^{3} \oplus \varepsilon^{0} \delta$ \\
\hline
\end{tabular}
\end{center}

\begin{center}
\begin{tabular}{|r|r|c|}
\hline
\multicolumn{3}{|l|}{\textbf{Level $N = 41$.  Nebentype $\eta = \chi_{41}^{20}$.  Field $\mathbb{F} = GF(21881^{60})$.}} \\
\multicolumn{3}{|l|}{Computed $T_{2}$, $T_{3}$, $T_{5}$.} \\
\hline
2 & $1$ & $\varepsilon^{0} \oplus \varepsilon^{1} \oplus \varepsilon^{2} \sigma_{41,2{f}}$ \\
\hline
2 & $1$ & $\varepsilon^{2} \oplus \varepsilon^{3} \oplus \varepsilon^{0} \sigma_{41,2{f}}$ \\
\hline
2 & $1$ & $\varepsilon^{0} \oplus \varepsilon^{1} \mathrm{Sym}^2(\sigma_{41,2{f}})$ \\
\hline
\end{tabular}
\end{center}

\subsection{} \label{resChi} For each~$N$,
the next table specifies the basis that Sage chooses for the group of
characters $(\Z/N\Z)^\times \to \F_p$.  If there is one basis element,
it is denoted~$\chi_N$.  If there is more than one, they are denoted
$\chi_{N,0}$, $\chi_{N,1}$, etc.  The \emph{order} of~$\chi$ is the
smallest positive~$n$ so that $\chi^n$ is trivial on
$(\Z/N\Z)^\times$.  The \emph{parity} is even if $\chi(-1) = +1$ and
odd if $\chi(-1) = -1$.

\begin{longtable}{|c|c|c|c|l|}
\hline
$\chi_{N,i}$ & $p$ & order & parity & definition \\
\hline
\hline
$\chi_{7}$ & 12037 & 6 & odd & $3 \mapsto -1293$  \\
\hline
$\chi_{9}$ & 12379 & 6 & odd & $2 \mapsto 5770$  \\
\hline
$\chi_{12,0}$ & 5413 & 2 & odd & $7 \mapsto -1$, $5 \mapsto 1$  \\
$\chi_{12,1}$ & 5413 & 2 & odd & $7 \mapsto 1$, $5 \mapsto -1$  \\
\hline
$\chi_{13}$ & 12037 & 12 & odd & $2 \mapsto 4019$  \\
\hline
$\chi_{15,0}$ & 12037 & 2 & odd & $11 \mapsto -1$, $7 \mapsto 1$  \\
$\chi_{15,1}$ & 12037 & 4 & odd & $11 \mapsto 1$, $7 \mapsto 3417$  \\
\hline
$\chi_{16,0}$ & 4001 & 2 & odd & $15 \mapsto -1$, $5 \mapsto 1$  \\
$\chi_{16,1}$ & 4001 & 4 & even & $15 \mapsto 1$, $5 \mapsto -899$  \\
\hline
$\chi_{17}$ & 16001 & 16 & odd & $3 \mapsto 83$  \\
\hline
$\chi_{18}$ & 3637 & 6 & odd & $11 \mapsto -695$  \\
\hline
$\chi_{19}$ & 3637 & 18 & odd & $2 \mapsto -31$  \\
\hline
$\chi_{20,0}$ & 12037 & 2 & odd & $11 \mapsto -1$, $17 \mapsto 1$  \\
$\chi_{20,1}$ & 12037 & 4 & odd & $11 \mapsto 1$, $17 \mapsto 3417$  \\
\hline
$\chi_{21,0}$ & 12037 & 2 & odd & $8 \mapsto -1$, $10 \mapsto 1$  \\
$\chi_{21,1}$ & 12037 & 6 & odd & $8 \mapsto 1$, $10 \mapsto -1293$  \\
\hline
$\chi_{22}$ & 16001 & 10 & odd & $13 \mapsto 3018$  \\
\hline
$\chi_{23}$ & 22067 & 22 & odd & $5 \mapsto 7863$  \\
\hline
$\chi_{24,0}$ & 12379 & 2 & odd & $7 \mapsto -1$, $13 \mapsto 1$, $17 \mapsto 1$  \\
$\chi_{24,1}$ & 12379 & 2 & even & $7 \mapsto 1$, $13 \mapsto -1$, $17 \mapsto 1$  \\
$\chi_{24,2}$ & 12379 & 2 & odd & $7 \mapsto 1$, $13 \mapsto 1$, $17 \mapsto -1$  \\
\hline
$\chi_{25}$ & 16001 & 20 & odd & $2 \mapsto 7734$  \\
\hline
$\chi_{26}$ & 12037 & 12 & odd & $15 \mapsto 4019$  \\
\hline
$\chi_{27}$ & 11863 & 18 & odd & $2 \mapsto 5034$  \\
\hline
$\chi_{28,0}$ & 12379 & 2 & odd & $15 \mapsto -1$, $17 \mapsto 1$  \\
$\chi_{28,1}$ & 12379 & 6 & odd & $15 \mapsto 1$, $17 \mapsto 5770$  \\
\hline
$\chi_{29}$ & 2297 & 28 & odd & $2 \mapsto 1108$  \\
\hline
$\chi_{31}$ & 4201 & 30 & odd & $3 \mapsto -1970$  \\
\hline
$\chi_{37}$ & 3889 & 36 & odd & $2 \mapsto -1338$  \\
\hline
$\chi_{41}$ & 21881 & 40 & odd & $6 \mapsto -10354$  \\
\hline
\end{longtable}

\subsection{} \label{resCusp}

In the following table we give the $q$-expansions of the holomorphic
cusp forms that we observed in our computations.  $S_{k}(N, \chi)$
denotes the space of weight~$k$ cusp forms on $\Gamma_{0}(N)$ with
character~$\chi$.  The notation~$\sigma_{N,k}$ for individual cusp
forms makes manifest the level~$N$ and weight~$k$.  The $q$-expansions
were computed using Sage \cite{sage}.

The field of definition of a cusp form, if not specified, is the field
generated by the coefficients we display. For instance, $q + 2iq^2 +
55q^3 + \cdots$ has coefficients in $\Q(i)$.  By~$\zeta_{m}$ we mean a
primitive $m$-th root of unity.  When we must specify the field, it is
in the line beginning with ``over''.

\begin{center}
\begin{longtable}{|l|}
\hline
$\sigma_{7,3} = q - 3q^{2} + 5q^{4} - 7q^{7} + O(q^{8})$ 
in $S_3(7, \chi_7^3)$ \\
\hline
$\sigma_{8,3} = q - 2q^{2} - 2q^{3} + 4q^{4} + 4q^{6} - 8q^{8} - 5q^{9} + 14q^{11} - 8q^{12} + O(q^{16})$ 
in $S_3(8, \chi_{24,0}\chi_{24,1})$ \\
\hline
$\sigma_{9,3} = q + \left(-\zeta_{6} - 1\right)q^{2} + \left(3 \zeta_{6} - 3\right)q^{3} - \zeta_{6}q^{4} + O(q^5)$ 
in $S_3(9, \chi_{27}^{15})$ \\
\hline
$\sigma_{11,2} = q - 2q^{2} - q^{3} + 2q^{4} + q^{5} + 2q^{6} - 2q^{7} + O(q^{8})$ 
in $S_2(11, 1)$ \\
\hline
$\sigma_{13,2} = q + \left(-\zeta_{6} - 1\right)q^{2} + \left(2 \zeta_{6} - 2\right)q^{3} + \zeta_{6}q^{4} + O(q^5)$ 
in $S_2(13, \chi_{13}^{2})$ \\
\hline
$\sigma_{13,4} = q - 5q^{2} - 7q^{3} + 17q^{4} - 7q^{5} + 35q^{6} - 13q^{7} + O(q^{8})$ 
in $S_4(13, 1)$ \\
\hline
$\sigma_{14,2} = q - q^{2} - 2q^{3} + q^{4} + 2q^{6} + q^{7} + O(q^{8})$ 
in $S_2(14, 1)$ \\
\hline
$\sigma_{15,2} = q - q^{2} - q^{3} - q^{4} + q^{5} + q^{6} + O(q^{8})$ 
in $S_2(15, 1)$ \\
\hline
$\sigma_{16,2} = q + \left(-i - 1\right)q^{2} + \left(i - 1\right)q^{3} + 2 iq^{4} + O(q^5)$ 
in $S_2(16, \chi_{16,1})$ \\
\hline
$\sigma_{17,2a} = q - q^{2} - q^{4} - 2q^{5} + 4q^{7} + O(q^{8})$ 
in $S_2(17, 1)$ \\
\hline
$\sigma_{17,2b} = q + \left(-\zeta_{8}^{3} + \zeta_{8}^{2} - 1\right)q^{2} + \left(\zeta_{8}^{3} - \zeta_{8}^{2} - \zeta_{8} - 1\right)q^{3} + O(q^4)$ 
in $S_2(17, \chi_{17}^{2})$ \\
\hline
$\sigma_{17,4} = q - 3q^{2} - 8q^{3} + q^{4} + 6q^{5} + 24q^{6} - 28q^{7} + O(q^{8})$ 
in $S_4(17, 1)$ \\
\hline
$\sigma_{18,2} = q - \zeta_{6}q^{2} + \left(\zeta_{6} - 2\right)q^{3} + \left(\zeta_{6} - 1\right)q^{4} + O(q^6)$ 
in $S_2(18, \chi_{18}^{2})$ \\
\hline
$\sigma_{19,2a} = q - 2q^{3} - 2q^{4} + 3q^{5} - q^{7} + O(q^{8})$ 
in $S_2(19, 1)$ \\
\hline
$\sigma_{19,2b} = q + \left(-\zeta_{18,2}^{2} + \zeta_{18,2} - 1\right)q^{2} + O(q^3)$ 
in $S_2(19, \chi_{19}^{2})$ \\
\hline
$\sigma_{19,4} = q - 3q^{2} - 5q^{3} + q^{4} - 12q^{5} + 15q^{6} + 11q^{7} + O(q^{8})$ 
in $S_4(19, 1)$ \\
\hline
$\sigma_{20,2a} = q - 2q^{3} - q^{5} + 2q^{7} + q^{9} + 2q^{13} + O(q^{14})$ 
in $S_2(20, 1)$ \\
\hline
$\sigma_{20,2b} = q + \left(-i - 1\right)q^{2} + 2 iq^{4} + \left(i - 2\right)q^{5} + O(q^8)$ 
in $S_2(20, \chi_{20,0}\chi_{20,1})$ \\
\hline
$\sigma_{21,2a} = q - q^{2} + q^{3} - q^{4} - 2q^{5} - q^{6} - q^{7} + O(q^{8})$ 
in $S_2(21, 1)$ \\
\hline
$\sigma_{21,2b} = q + \left(2 \zeta_{6} - 2\right)q^{2} - \zeta_{6}q^{3} - 2 \zeta_{6}q^{4} + \left(-2 \zeta_{6} + 2\right)q^{5} + O(q^6)$ 
in $S_2(21, \chi_{21,1}^{2})$ \\
\hline
$\sigma_{21,2c} = q + \left(-\zeta_{6} - 1\right)q^{3} + \left(2 \zeta_{6} - 2\right)q^{4} + O(q^{7})$
in $S_2(21, \chi_{21,0}\chi_{21,1})$ \\
\hline
$\sigma_{21,4} = q - 3q^{2} - 3q^{3} + q^{4} - 18q^{5} + 9q^{6} + 7q^{7} + O(q^{8})$ 
in $S_4(21, 1)$ \\
\hline
$\sigma_{22,2} = q - \zeta_{10}q^{2} + \left(-\zeta_{10}^{3} + \zeta_{10} - 1\right)q^{3} + \zeta_{10}^{2}q^{4} + O(q^5)$ 
in $S_2(22, \chi_{22}^{2})$ \\
\hline
$\sigma_{22,4} = q - 2q^{2} - 7q^{3} + 4q^{4} - 19q^{5} + 14q^{6} + 14q^{7} + O(q^{8})$ 
in $S_4(22, 1)$ \\
\hline
$\sigma_{23,2a} = q + b_{0}q^{2} + \left(-2 b_{0} - 1\right)q^{3} + \left(-b_{0} - 1\right)q^{4} + 2 b_{0}q^{5} + O(q^6)$ 
in $S_2(23, 1)$ \\
$\mathrm{over\ }\mathbb{Q}[b_{0}]/(b_{0}^{2} + b_{0} - 1)$ \\
\hline
$\sigma_{23,2b} = q + \left(\zeta_{22}^{9} - \zeta_{22}^{6} - \zeta_{22}^{4} - 1\right)q^{2} + O(q^3)$
in $S_2(23, \chi_{23}^{2})$ \\
\hline
$\sigma_{23,4} = q - 2q^{2} - 5q^{3} - 4q^{4} - 6q^{5} + 10q^{6} - 8q^{7} + O(q^{8})$ 
in $S_4(23, 1)$ \\
\hline
$\sigma_{24,2a} = q - q^{3} - 2q^{5} + q^{9} + 4q^{11} - 2q^{13} + O(q^{14})$ 
in $S_2(24, 1)$ \\
\hline
$\sigma_{24,2b} = q + b_{0}q^{2} + \left(b_{0} + 1\right)q^{3} + \left(-2 b_{0} - 2\right)q^{4} + O(q^5)$ 
in $S_2(24, \chi_{24,1})$ \\
$\mathrm{over\ }\mathbb{Q}[b_{0}]/(b_{0}^{2} + 2 b_{0} + 2)$ \\
\hline
$\sigma_{24,2c} = q + b_{0}q^{2} + \left(-b_{0} - 1\right)q^{3} - 2q^{4} + \left(-b_{0} + 2\right)q^{6} + O(q^8)$ 
in $S_2(24, \chi_{24,0}\chi_{24,1}\chi_{24,2})$ \\
$\mathrm{over\ }\mathbb{Q}[b_{0}]/(b_{0}^{2} + 2)$ \\
\hline
$\sigma_{25,2a} = q + b_{0}q^{2} + \left(\left(\zeta_{10}^{3} + \zeta_{10} - 1\right) b_{0} + \zeta_{10}^{2} - 1\right)q^{3} + O(q^4)$ 
in $S_2(25, \chi_{25}^{2})$ \\
$\mathrm{over\ }( \mathbb{Q}(\zeta_{10}) )[b_{0}]/(b_{0}^{2} + \left(\zeta_{10} + 1\right) b_{0} + \zeta_{10}^{2} - 2 \zeta_{10} + 1)$ \\
\hline
$\sigma_{25,2b} = q + \left(-\zeta_{5}^{3} - \zeta_{5} - 1\right)q^{2} + \zeta_{5}q^{3} + \left(-\zeta_{5}^{2} - \zeta_{5} - 1\right)q^{4} + O(q^5)$ 
in $S_2(25, \chi_{25}^{4})$ \\
\hline
$\sigma_{25,4} = q - q^{2} - 7q^{3} - 7q^{4} + 7q^{6} - 6q^{7} + O(q^{8})$ 
in $S_4(25, 1)$ \\
\hline
$\sigma_{26,2a} = q - q^{2} + q^{3} + q^{4} - 3q^{5} - q^{6} - q^{7} + O(q^{8})$ 
in $S_2(26, 1)$ \\
\hline
$\sigma_{26,2b} = q + q^{2} - 3q^{3} + q^{4} - q^{5} - 3q^{6} + q^{7} + O(q^{8})$ 
in $S_2(26, 1)$ \\
\hline
$\sigma_{26,2c} = q + \left(-\zeta_{3} - 1\right)q^{2} + \zeta_{3}q^{4} - q^{5} + 4 \zeta_{3}q^{7} + O(q^{8})$ 
in $S_2(26, \chi_{26}^{4})$ \\
\hline
$\sigma_{26,2d} = q + b_{0}q^{2} - q^{3} - q^{4} - 3 b_{0}q^{5} - b_{0}q^{6} + 3 b_{0}q^{7} + O(q^{8})$ 
in $S_2(26, \chi_{26}^{6})$ \\
$\mathrm{over\ }\mathbb{Q}[b_{0}]/(b_{0}^{2} + 1)$ \\
\hline
$\sigma_{27,2a} = q - 2q^{4} - q^{7} + O(q^{8})$ 
in $S_2(27, 1)$ \\
\hline
$\sigma_{27,2b} = q + b_{0}q^{2} + \left(\left(\zeta_{18,2}^{5} - \zeta_{18,2}\right) b_{0} - \zeta_{18,2}^{3} + \zeta_{18,2}^{2} - \zeta_{18,2}\right)q^{3} + O(q^4)$ 
in $S_2(27, \chi_{27}^{2})$ \\
$\mathrm{over\ }( \mathbb{Q}(\zeta_{18,2}) )[b_{0}]/(b_{0}^{2} + \left(\zeta_{18,2}^{2} - \zeta_{18,2} + 1\right) b_{0} + \zeta_{18,2}^{4} - \zeta_{18,2}^{3} - \zeta_{18,2}^{2} - \zeta_{18,2} + 1)$ \\
\hline
$\sigma_{27,4} = q - 3q^{2} + q^{4} - 15q^{5} - 25q^{7} + O(q^{8})$ 
in $S_4(27, 1)$ \\
\hline
$\sigma_{28,2a} = q - \zeta_{6}q^{3} + \left(3 \zeta_{6} - 3\right)q^{5} + \left(-2 \zeta_{6} - 1\right)q^{7} + O(q^{8})$ 
in $S_2(28, \chi_{28,1}^{2})$ \\
\hline
$\sigma_{28,2b} = q + b_{0}q^{2} + \left(\left(\zeta_{6} - 2\right) b_{0} - \zeta_{6} - 1\right)q^{3} + O(q^4)$ 
in $S_2(28, \chi_{28,0}\chi_{28,1})$ \\
$\mathrm{over\ }( \mathbb{Q}(\zeta_{6}) )[b_{0}]/(b_{0}^{2} + 2 \zeta_{6} b_{0} + 2 \zeta_{6} - 2)$ \\
\hline
$\sigma_{28,2c} = q + b_{0}q^{2} + \left(-b_{0} - 2\right)q^{4} + \left(-2 b_{0} - 1\right)q^{7} + O(q^{8})$ 
in $S_2(28, \chi_{28,0}\chi_{28,1}^{3})$ \\
$\mathrm{over\ }\mathbb{Q}[b_{0}]/(b_{0}^{2} + b_{0} + 2)$ \\
\hline
$\sigma_{28,4} = q - 10q^{3} - 8q^{5} - 7q^{7} + 73q^{9} - 40q^{11} - 12q^{13} + O(q^{14})$ 
in $S_4(28, 1)$ \\
\hline
$\sigma_{29,2a} = q + b_{0}q^{2} - b_{0}q^{3} + \left(-2 b_{0} - 1\right)q^{4} - q^{5} + \left(2 b_{0} - 1\right)q^{6} + O(q^7)$ 
in $S_2(29, 1)$ \\
$\mathrm{over\ }\mathbb{Q}[b_{0}]/(b_{0}^{2} + 2 b_{0} - 1)$ \\
\hline
$\sigma_{29,2b} = q + b_{0}q^{2} + \left(\left(\zeta_{14}^{3} + \zeta_{14}^{2} + \zeta_{14}\right) b_{0} + \zeta_{14}^{4} + \zeta_{14}^{3} + \zeta_{14}^{2} + \zeta_{14}\right)q^{3} + O(q^4)$ 
in $S_2(29, \chi_{29}^{2})$ \\
$\mathrm{over\ }( \mathbb{Q}(\zeta_{14}) )[b_{0}]/(b_{0}^{2} + \left(-\zeta_{14}^{5} + \zeta_{14}^{3} + \zeta_{14} + 1\right) b_{0} - \zeta_{14}^{5} + \zeta_{14}^{4} + \zeta_{14}^{2} - \zeta_{14} + 1)$ \\
\hline
$\sigma_{29,2c} = q + \left(-\zeta_{7}^{5} - \zeta_{7}^{4} - \zeta_{7}^{3} - \zeta_{7} - 1\right)q^{2} + \left(-\zeta_{7}^{5} - 1\right)q^{3} + O(q^4)$ 
in $S_2(29, \chi_{29}^{4})$ \\
\hline
$\sigma_{29,2d} = q + b_{0}q^{2} - b_{0}q^{3} - 3q^{4} - 3q^{5} + 5q^{6} + 2q^{7} + O(q^{8})$ 
in $S_2(29, \chi_{29}^{14})$ \\
$\mathrm{over\ }\mathbb{Q}[b_{0}]/(b_{0}^{2} + 5)$ \\
\hline
$\sigma_{29,4} = q + b_{0}q^{2} + \left(-3 b_{0} - 8\right)q^{3} + \left(-2 b_{0} - 7\right)q^{4} + \left(4 b_{0} - 1\right)q^{5} + O(q^6)$ 
in $S_4(29, 1)$ \\
$\mathrm{over\ }\mathbb{Q}[b_{0}]/(b_{0}^{2} + 2 b_{0} - 1)$ \\
\hline
$\sigma_{31,2a} = q + b_{0}q^{2} - 2 b_{0}q^{3} + \left(b_{0} - 1\right)q^{4} + q^{5} + \left(-2 b_{0} - 2\right)q^{6} + O(q^{7})$ 
in $S_2(31, 1)$ \\
$\mathrm{over\ }\mathbb{Q}[b_{0}]/(b_{0}^{2} - b_{0} - 1)$ \\
\hline
$\sigma_{31,2b} = q + b_{0}q^{2} + \left(\left(\zeta_{30}^{5} - 2 \zeta_{30}^{3} - \zeta_{30}^{2} + \zeta_{30} + 1\right) b_{0} + \zeta_{30}^{6} + \zeta_{30}^{5} - \zeta_{30}^{4} - 2 \zeta_{30}^{3} - \zeta_{30}^{2} + \zeta_{30}\right)q^{3}$ \\
$ + \,O(q^4)$ 
in $S_2(31, \chi_{31}^{2})$ \\
$\mathrm{over\ }( \mathbb{Q}(\zeta_{30}) )[b_{0}]/(b_{0}^{2} + \left(-\zeta_{30}^{3} + 1\right) b_{0} + 2 \zeta_{30}^{6} - \zeta_{30}^{4} + \zeta_{30}^{3} - \zeta_{30}^{2} + 2)$ \\
\hline
$\sigma_{31,2c} = q + \left(\zeta_{5}^{3} + \zeta_{5}^{2} + \zeta_{5}\right)q^{2} - \zeta_{5}^{3}q^{3} + \left(\zeta_{5}^{3} + 1\right)q^{4} + O(q^5)$ 
in $S_2(31, \chi_{31}^{6})$ \\
\hline
$\sigma_{31,2d} = q + b_{0}q^{2} + \left(\left(-\zeta_{3} - 1\right) b_{0}\right)q^{3} + \left(-2 b_{0} - 1\right)q^{4} + O(q^5)$ 
in $S_2(31, \chi_{31}^{10})$ \\
$\mathrm{over\ }( \mathbb{Q}(\zeta_{3}) )[b_{0}]/(b_{0}^{2} + 2 b_{0} - 1)$ \\
\hline
$\sigma_{31,4} = q + b_{0}q^{2} + \left(-2 b_{0} - 6\right)q^{3} + \left(-5 b_{0} - 10\right)q^{4} + \left(3 b_{0} - 5\right)q^{5} + O(q^{6})$ 
in $S_4(31, 1)$ \\
$\mathrm{over\ }\mathbb{Q}[b_{0}]/(b_{0}^{2} + 5 b_{0} + 2)$ \\
\hline
$\sigma_{37,2a} = q - 2q^{2} - 3q^{3} + 2q^{4} - 2q^{5} + 6q^{6} - q^{7} + 6q^{9} + 4q^{10} - 5q^{11} + O(q^{12})$ 
in $S_2(37, 1)$ \\
\hline
$\sigma_{37,2b} = q + q^{3} - 2q^{4} - q^{7} - 2q^{9} + 3q^{11} - 2q^{12} - 4q^{13} + O(q^{14})$ 
in $S_2(37, 1)$ \\
\hline
$\sigma_{37,2c} = q + b_{0}q^{2}$ \\
$ + \left(\left(\zeta_{18,2}^{4} + \zeta_{18,2}^{2}\right) b_{0}^{2} + \left(\zeta_{18,2}^{5} + \zeta_{18,2}^{4} + \zeta_{18,2}^{3} + \zeta_{18,2}^{2}\right) b_{0} - 2 \zeta_{18,2}^{5} + \zeta_{18,2}^{4} - \zeta_{18,2}^{3} + \zeta_{18,2}^{2} - 1\right)q^{3}$ \\
$+\,O(q^4)$ 
in $S_2(37, \chi_{37}^{2})$ \\
$\mathrm{over\ }( \mathbb{Q}(\zeta_{18,2}) )[b_{0}]/(b_{0}^{3} + \left(-\zeta_{18,2}^{4} + \zeta_{18,2}^{3} + 2 \zeta_{18,2} + 1\right) b_{0}^{2}$ \\
$\qquad + \left(-2 \zeta_{18,2}^{5} + 2 \zeta_{18,2}^{3} + 2 \zeta_{18,2}^{2} - 2 \zeta_{18,2}\right) b_{0} + \zeta_{18,2}^{5} - \zeta_{18,2}^{4} + \zeta_{18,2}^{3} - 2 \zeta_{18,2}^{2} - \zeta_{18,2} + 1)$ \\
\hline
$\sigma_{37,2d} = q + \left(-\zeta_{9} - 1\right)q^{2} + \left(-\zeta_{9}^{4} + \zeta_{9}^{3} - \zeta_{9}^{2} + 1\right)q^{3} + O(q^4)$ 
in $S_2(37, \chi_{37}^{4})$ \\
\hline
$\sigma_{37,2e} = q + \left(-\zeta_{9}^{5} + \zeta_{9}^{4} - \zeta_{9}^{3} + \zeta_{9}\right)q^{2} + \left(\zeta_{9}^{5} + \zeta_{9}^{2} - 1\right)q^{3} + O(q^4)$ 
in $S_2(37, \chi_{37}^{4})$ \\
\hline
$\sigma_{37,2{f}} = q + b_{0}q^{2} + \left(\left(-\zeta_{6} - 1\right) b_{0} + \zeta_{6}\right)q^{3} - \zeta_{6}q^{4} + O(q^5)$ 
in $S_2(37, \chi_{37}^{6})$ \\
$\mathrm{over\ }( \mathbb{Q}(\zeta_{6}) )[b_{0}]/(b_{0}^{2} - \zeta_{6})$ \\
\hline
$\sigma_{37,2g} = q + \left(-\zeta_{3} - 1\right)q^{2} - \zeta_{3}q^{4} + \zeta_{3}q^{5} + 2 \zeta_{3}q^{7} + O(q^{8})$ 
in $S_2(37, \chi_{37}^{12})$ \\
\hline
$\sigma_{37,2h} = q + b_{0}q^{2} - q^{3} - 2q^{4} - b_{0}q^{5} - b_{0}q^{6} + 3q^{7} + O(q^{8})$ 
in $S_2(37, \chi_{37}^{18,2})$ \\
$\mathrm{over\ }\mathbb{Q}[b_{0}]/(b_{0}^{2} + 4)$ \\
\hline
$\sigma_{37,4} = q + b_{0}q^{2} + \left(-\frac{1}{8} b_{0}^{3} - \frac{9}{8} b_{0}^{2} - \frac{13}{4} b_{0} - \frac{11}{4}\right)q^{3} + \left(b_{0}^{2} - 8\right)q^{4} + O(q^5)$ 
in $S_4(37, 1)$ \\
$\mathrm{over\ }\mathbb{Q}[b_{0}]/(b_{0}^{4} + 6 b_{0}^{3} - b_{0}^{2} - 16 b_{0} + 6)$ \\
\hline
$\sigma_{41,2a} = q + b_{0}q^{2} + \left(-\frac{1}{2} b_{0}^{2} - b_{0} + \frac{3}{2}\right)q^{3} + \left(b_{0}^{2} - 2\right)q^{4} + O(q^5)$ 
in $S_2(41, 1)$ \\
$\mathrm{over\ }\mathbb{Q}[b_{0}]/(b_{0}^{3} + b_{0}^{2} - 5 b_{0} - 1)$ \\
\hline
$\sigma_{41,2b} = q + b_{0}q^{2}$ \\
$ + \left(\left(\zeta_{20}^{3} - \zeta_{20}^{2} - \zeta_{20} + 1\right) b_{0}^{2} + \left(-\zeta_{20}^{7} + \zeta_{20}^{6} + \zeta_{20}^{3} - \zeta_{20}^{2}\right) b_{0} \right.$ \\
$\qquad\left. + 2 \zeta_{20}^{5} - 2 \zeta_{20}^{4} - \zeta_{20}^{3} + \zeta_{20}^{2} + 2 \zeta_{20} - 2\right)q^{3} + O(q^4)$ 
in $S_2(41, \chi_{41}^{18,2})$ \\
$\mathrm{over\ } ( \mathbb{Q}(\zeta_{20}) )[b_{0}]/(b_{0}^{3} + \left(\zeta_{20}^{6} + \zeta_{20}^{3} + 1\right) b_{0}^{2}$ \\
$\qquad + \left(\zeta_{20}^{7} - 3 \zeta_{20}^{6} - \zeta_{20}^{5} + 2 \zeta_{20}^{3} - \zeta_{20}^{2} - \zeta_{20} + 1\right) b_{0} - 2 \zeta_{20}^{6} + \zeta_{20}^{5} + \zeta_{20}^{4} - \zeta_{20}^{3} + 2)$ \\
\hline
$\sigma_{41,2c} = q + b_{0}q^{2} + \left(\left(\frac{2}{5} \zeta_{10}^{3} + \frac{1}{5} \zeta_{10}^{2} - \frac{4}{5} \zeta_{10} + \frac{2}{5}\right) b_{0} + \frac{2}{5} \zeta_{10}^{3} + \frac{6}{5} \zeta_{10}^{2} - \frac{4}{5} \zeta_{10} + \frac{2}{5}\right)q^{3}$ \\
$ + \,O(q^4)$ 
in $S_2(41, \chi_{41}^{4})$ \\
$\mathrm{over\ }( \mathbb{Q}(\zeta_{10}) )[b_{0}]/(b_{0}^{2} + \left(-\zeta_{10} + 1\right) b_{0} + \zeta_{10}^{2} + \zeta_{10} + 1)$ \\
\hline
$\sigma_{41,2d} = q + b_{0}q^{2} + \zeta_{5}^{2} b_{0}q^{3} + \left(\left(-2 \zeta_{5}^{3} - \zeta_{5} - 1\right) b_{0} + \zeta_{5}^{2} - \zeta_{5} + 1\right)q^{4} + O(q^5)$ 
in $S_2(41, \chi_{41}^{8})$ \\
$\mathrm{over\ }( \mathbb{Q}(\zeta_{5}) )[b_{0}]/(b_{0}^{2} + \left(2 \zeta_{5}^{3} + \zeta_{5} + 1\right) b_{0} - \zeta_{5}^{2} - \zeta_{5} - 1)$ \\
\hline
$\sigma_{41,2e} = q + b_{0}q^{2} + \left(\left(\frac{1}{2} i - \frac{1}{2}\right) b_{0}^{2} + \frac{5}{2} i - \frac{5}{2}\right)q^{3} + O(q^4)$ 
in $S_2(41, \chi_{41}^{10})$ \\
$\mathrm{over\ }( \mathbb{Q}(i) )[b_{0}]/(b_{0}^{3} - i b_{0}^{2} + 5 b_{0} - 3 i)$ \\
\hline
$\sigma_{41,2{f}} = q - q^{2} + \left(-\frac{1}{2} b_{0} - \frac{1}{2}\right)q^{3} - q^{4} + 2q^{5} + \left(\frac{1}{2} b_{0} + \frac{1}{2}\right)q^{6} +O(q^7)$ 
in $S_2(41, \chi_{41}^{20})$ \\
$\mathrm{over\ }\mathbb{Q}[b_{0}]/(b_{0}^{2} + 2 b_{0} + 33)$ \\
\hline
$\sigma_{41,4} = q + b_{0}q^{2} + \left(-\frac{1}{2} b_{0}^{2} - 3 b_{0} - \frac{5}{2}\right)q^{3} + \left(b_{0}^{2} - 8\right)q^{4} + O(q^5)$ 
in $S_4(41, 1)$ \\
$\mathrm{over\ }\mathbb{Q}[b_{0}]/(b_{0}^{3} + 3 b_{0}^{2} - 5 b_{0} - 3)$ \\
\hline
\end{longtable}
\end{center}

%%%%%%%%%%%%%%%%%%%%%%%%%%%%%%%%%%%%%%%%%%%%%%%%%%%%%%%%%%%%

\bibliography{AGM-VII}

\end{document}